\documentclass[12pt]{amsart}

\usepackage{times,epsf,amssymb,amsmath,hyperref, yhmath, rlepsf, color}

\begin{document}

\newtheorem{thm}{Theorem}[section]
\newtheorem{lem}[thm]{Lemma}
\newtheorem{cor}[thm]{Corollary}
\newtheorem{conj}[thm]{Conjecture}

\theoremstyle{definition}
\newtheorem{defn}[thm]{\bf{Definition}}

\theoremstyle{remark}
\newtheorem{rmk}[thm]{Remark}

\theoremstyle{question}
\newtheorem{question}[thm]{Question}

\def\square{\hfill${\vcenter{\vbox{\hrule height.4pt \hbox{\vrule width.4pt height7pt \kern7pt \vrule width.4pt} \hrule height.4pt}}}$}

\newcommand{\SI}{\partial_\infty ({\Bbb H}^2\times {\Bbb R})}
\newcommand{\Si}{S^1_{\infty}\times {\Bbb R}}
\newcommand{\si}{S^1_{\infty}({\Bbb H}^2)}
\newcommand{\CS}{S^1_\infty\times\{\pm\infty\}}
\newcommand{\PI}{\partial_{\infty}}
\newcommand{\caps}{{\Bbb H}^2\times\{\pm\infty\}}

\newcommand{\BH}{\Bbb H}
\newcommand{\BHH}{{\Bbb H}^2\times {\Bbb R}}
\newcommand{\BR}{\Bbb R}
\newcommand{\BC}{\Bbb C}
\newcommand{\BZ}{\Bbb Z}

\newcommand{\wh}{\widehat}
\newcommand{\wt}{\widetilde}
\newcommand{\wc}{\check}

\newcommand{\A}{\mathcal{A}}
\newcommand{\Oo}{\mathcal{O}}
\newcommand{\V}{\mathcal{V}}
\newcommand{\U}{\mathcal{U}}
\newcommand{\Y}{\mathcal{Y}}
\newcommand{\W}{\mathcal{W}}
\newcommand{\C}{\mathcal{C}}
\newcommand{\K}{\mathcal{K}}
\newcommand{\D}{\mathcal{D}}
\newcommand{\p}{\mathcal{P}}
\newcommand{\R}{\mathcal{R}}
\newcommand{\h}{\mathcal{H}}
\newcommand{\T}{\mathcal{T}}
\newcommand{\N}{\mathcal{N}}
\newcommand{\s}{\mathcal{S}}

\newcommand{\e}{\epsilon}
\newcommand{\B}{\mathbf{B}}

\newenvironment{pf}{{\it Proof:}\quad}{\square \vskip 12pt}

\title[Asymptotic Plateau Problem in $\BHH$]{Asymptotic Plateau Problem in $\BHH$}
\author{Baris Coskunuzer}
\address{Boston College, Mathematics Department, Chestnut Hill, MA 02467}
\email{coskunuz@bc.edu}
\thanks{The author is partially supported by BAGEP award of the Science Academy, and a Royal Society Newton Mobility Grant.}

\maketitle

\begin{abstract}

We give a fairly complete solution to the asymptotic Plateau Problem for area minimizing surfaces in $\BHH$. In particular, we identify the collection of Jordan curves in $\SI$ which bounds an area minimizing surface in $\BHH$. Furthermore, we study the similar problem for minimal surfaces, and show that the situation is highly different.
\end{abstract}

\section{Introduction}

Asymptotic Plateau Problem asks the existence of a minimal surface $\Sigma$ in $\BHH$ for a given curve $\Gamma$ in $\SI$ with $\PI\Sigma=\Gamma$.
We will call a finite collection of disjoint Jordan curves $\Gamma$ in $\SI$ {\em fillable}, if $\Gamma$ bounds a complete, embedded \underline{minimal} surface $S$ in $\BHH$ with $\PI S=\Gamma$. We will call $\Gamma$ {\em strongly fillable} if $\Gamma$ bounds a complete, embedded, \underline{area minimizing} surface $\Sigma$ in $\BHH$ with $\PI \Sigma=\Gamma$. In this paper, our aim is to classify fillable and strongly fillable infinite curves in $\SI$.

In the last decade, minimal surfaces in $\BHH$ have been studied extensively, and many important results have been  obtained on the existence of many different types of minimal surfaces in $\BHH$ and their properties, e.g. \cite{NR, CR, CMT, MMR, MoR, MRR, PR, RT, ST, ST2}.

Let $\Gamma^\pm=\Gamma\cap (\overline{\BH^2}\times\{\pm\infty\}$) and $\wt{\Gamma}=\Gamma-(\Gamma^+\cup\Gamma^-)$. We call $\Gamma$ {\em infinite} if either $\Gamma^+$ or $\Gamma^-$ is nonempty, and {\em finite} otherwise. We will call $\Gamma$ {\em tame}, if $\Gamma^\pm$ consist of finitely many components.

Recently, Kloeckner and Mazzeo \cite{KM}, and the author \cite{Co} independently studied the asymptotic Plateau problem in $\BHH$.
While the author gave a fairly complete classification of strongly fillable finite curves in \cite{Co}, Kloeckner and Mazzeo constructed many interesting families of finite and infinite fillable curves in \cite{KM}. In this paper, we study the infinite curves in $\SI$, and complete the classification of strongly fillable curves. Our main result is as follows:

\begin{thm} Let $\Gamma$ be a tame infinite curve in $\SI$. Then, $\Gamma$ is strongly fillable if and only if all of the following conditions satisfied:
	\begin{itemize}
		\item $int(\Gamma^\pm)$ is a collection of geodesics (possibly empty).
		\item $\Gamma$ is tall.	
		\item $\Gamma$ is nonoverlapping at the corner.
		\item $\Gamma$ is fat at infinity.
	\end{itemize}
\end{thm}

 
For a given $\Gamma$, being {\em fat at infinity} is completely determined by $\Gamma^\pm$, while being {\em tall} is determined by $\wt{\Gamma}$. In \cite{Co}, we introduced the notion of being tall, and showed that a finite curve $\Gamma$ is strongly fillable if and only if $\Gamma$ is tall. In this paper, we introduce two new notions for infinite curves in $\SI$, namely being \textit{fat or skinny at infinity}, and {\em nonoverlapping at the corner}, which characterize the behavior of $\Gamma^\pm$. Note that we have some exceptional curves coming from fat at infinity definition (See Remark \ref{exception}).

The outline of the paper is as follows. The main difference for the solution to the asymptotic Plateau problem for $\BH^3$ and for $\BHH$ is the {\em escaping to infinity} problem. In $\BH^3$, a sequence of compact area minimizing surfaces $\{\Sigma_n\}$ with $\partial \Sigma_n\to\Gamma$ limits to an area minimizing surface $\Sigma$ with $\PI\Sigma=\Gamma$ since convex hull of $\Gamma$ is a natural barrier for the sequence $\{\Sigma_n\}$ \cite{An}. However, in $\BHH$, such a sequence might escape to the infinity, and give an empty limit as discussed in \cite{Co}. Hence, for the existence part of the main result, we first build a barrier near infinity to cover the complement of $\Gamma$ in $\SI$ so that no piece of the sequence $\Sigma_n$ escapes to infinity, and the limit area minimizing surface $\Sigma$ has the desired asymptotic boundary, i.e. $\PI \Sigma=\Gamma$. To build the barrier near infinity, we use finite tall rectangles for the cylinder $\Si$, and infinite rectangles and Scherk graphs (see section \ref{scherksec}) for the caps $\overline{\BH^2}\times\{\pm\infty\}$. For the nonexistence direction, we use the Scherk graphs and area comparison to get a contradiction.

In the proof of the {\em only if} part of the main result, the area comparison is crucial, so the result only applies to area minimizing surfaces in $\BHH$. In Section 4, we discuss many infinite Jordan curves which bounds complete minimal surfaces, but no area minimizing surfaces in $\BHH$.

The organization of the paper is as follows. In the next section, we give some definitions and related results. In Section 3, we give the complete solution to the asymptotic Plateau problem in $\BHH$, and prove the main result, Theorem 1.1. In Section 4, we discuss fillable and nonfillable infinite curves in $\BHH$. Finally in section 5, we give some concluding remarks on further generalizations and directions.

\subsection{Acknowledgements}

Part of this research was carried out at MIT during my visit. I would like to thank them for their great hospitality. I would like to thank the referee for very valuable remarks.

\section{Preliminaries}

In this section, we will give the basic definitions, and a brief overview of the past results, which will be used in the paper. 

Throughout the paper, we use the product compactification of $\BHH$. In particular, $\overline{\BHH}=\overline{\BH^2}\times\overline{\BR}=\BHH\cup \SI$ where $\SI$ consists of three components, i.e. the infinite cylinder $\Si$ and the caps at infinity $\overline{\BH^2}\times\{+\infty\}$, $\overline{\BH^2}\times\{-\infty\}$. Hence, $\overline{\BHH}$ is a solid cylinder under this compactification.


Let $\Sigma$ be an open, complete surface in $\BHH$, and $\PI\Sigma$ represent the asymptotic boundary of $\Sigma$ in $\SI$. Then, if $\overline{\Sigma}$ is the closure of $\Sigma$ in $\overline{\BHH}$, then $\PI \Sigma= \overline{\Sigma}\cap \SI$.

\begin{defn} A surface is {\em minimal} if the mean curvature $H$ vanishes everywhere. A compact surface with boundary $\Sigma$ is called {\em area minimizing surface} if $\Sigma$ has the smallest area among the surfaces with the same boundary. A noncompact surface is called {\em area minimizing surface} if any compact subsurface is an area minimizing surface. 
\end{defn}


Note that any area minimizing surface is minimal. In this paper, we will study the Jordan curves in $\SI$ which bounds a complete, embedded, minimal or area minimizing surfaces in $\BHH$.

\vspace{.2cm}

\noindent{\bf Convention:} [Curve] By \textit{curve}, we mean a finite collection of disjoint Jordan curves in $\SI$ throughout the paper unless otherwise stated.

\begin{defn} [Fillable Curves] Let $\Gamma$ be a curve in $\SI$. We will call $\Gamma$ {\em fillable} if $\Gamma$ bounds a complete embedded minimal surface $S$ in $\BHH$, i.e. $\PI S=\Gamma$. We will call $\Gamma$ {\em strongly fillable} if $\Gamma$ bounds a complete embedded area minimizing surface $\Sigma$ in $\BHH$, i.e. $\PI \Sigma=\Gamma$. We call such $S$ or $\Sigma$ as {\em filling surface} for $\Gamma$. 
\end{defn}

Notice that a strongly fillable curve is fillable. Note also that fillable curves here corresponds to {\em minimally fillable curves} in \cite{KM}.

The Asymptotic Plateau Problem for $\BHH$ is the following classification problems:

\vspace{.2cm}

{\em Which $\Gamma$ in $\SI$ is fillable or strongly fillable?}

\vspace{.2cm}

Throughout the paper, we will use the following notation for the curves at infinity. $\Gamma=\Gamma^+\cup\Gamma^-\cup\wt{\Gamma}$ where
$\Gamma^\pm=\Gamma\cap (\overline{\BH^2}\times\{\pm\infty\})$ and $\wt{\Gamma}=\Gamma\cap(\Si)$. In particular, $\Gamma^\pm$ is a collection of closed arcs and points in the caps at infinity, where $\wt{\Gamma}$ is a collection of open arcs and closed curves in the infinite cylinder.

With this notation, we will call a curve $\Gamma$ {\em finite} if $\Gamma^+=\Gamma^-=\emptyset$. We will call $\Gamma$ {\em infinite} otherwise.

\subsection{Finite Curves} \

When $\Gamma$ is an essential Jordan curve in $\Si$ which is a vertical graph over $S^1_\infty\times\{0\}$, then there exists a vertical graph over $\BH^2\times\{0\}$ giving a positive answer to this existence question \cite{NR}. However, for some nullhomotopic simple closed curves in $\SI$, the situation can be quite different. Unlike the $\BH^3$ case \cite{An}, Sa Earp and Toubiana  proved that there are some nonfillable $\Gamma$ in $\SI$ \cite{ST}.

\begin{defn} \label{tail} [Thin tail] Let $\Gamma$ be a simple closed curve in $\SI$, and let $\gamma$ be an arc in $\Gamma$. Assume that there is a vertical straight line $L_0$ in $\SI$ such that

\begin{itemize}

\item $\gamma \cap L_0 \neq \emptyset$ and $\partial \gamma\cap L_0 = \emptyset$,

\item $\gamma$ stays in one side of $L_0$,

\item $\gamma\subset \partial_\infty \BH^2 \times (c,c+\pi)$ for some $c\in \BR$.

\end{itemize}

Then, we call $\gamma$ {\em a thin tail} in $\Gamma$.

\end{defn}

\begin{lem} \cite{ST} \label{thin1} Let $\Gamma$ be a simple closed curve in $\SI$. If $\Gamma$ contains a thin tail, then there is no properly immersed minimal surface $\Sigma$ in $\BHH$ with $\PI \Sigma\supset\Gamma$.
\end{lem}


The above result shows that the curves with thin tail cannot be fillable. Hence, to bypass this obstruction, we introduced the following notion.

\begin{defn} \label{talldef} [Tall Curves] \cite{Co} Consider $\Si$ with the coordinates $(\theta, t)$ where $\theta\in [0,2\pi)$ and $t\in \BR$. We will call the rectangle $R=[\theta_1,\theta_2]\times[t_1,t_2]\subset \Si$ as {\em tall rectangle} if $t_2-t_1> \pi$.

We call a curve $\Gamma$ in $\SI$  {\em tall curve} if the open region $\Si-\Gamma$ can be written as a union of tall rectangles $int(R_i)$, i.e. $\Si-\Gamma=\bigcup_i int(R_i)$. Note that this definition naturally generalizes to infinite curves.
\end{defn}

Notice that tall curves do not have thin tails. Furthermore, we define the {\em height of a curve}, $h(\Gamma)$, as the length of the smallest component in vertical line segments in $\Si-\Gamma$. Hence, $\Gamma$ is tall if and only if $h(\Gamma)>\pi$.




In \cite{Co}, we gave a fairly complete classification of strongly fillable finite curves as follows.

\begin{lem} \cite{Co} \label{finiteAPP} Let $\Gamma$ be a finite curve in $\Si$ with $h(\Gamma)\neq \pi$. Then, $\Gamma$ is strongly fillable if and only if $\Gamma$ is a tall curve.
\end{lem}


After this result for finite curves, we aim to give a characterization for \textit{infinite} strongly fillable curves to complete the classification.

\subsection{Infinite Curves} \

As we defined before if $\Gamma^\pm=\Gamma\cap (\overline{\BH^2}\times\{\pm\infty\}$) is nonempty, we call $\Gamma$ in $\SI$ an infinite curve. In \cite{KM}, Kloeckner and Mazzeo studied infinite fillable curves. To cite their result, we need to adapt their notation.

Note that here we use closed caps $\overline{\BH^2}\times\{\pm\infty\}$ to define infinite curves. However, Kloeckner-Mazzeo used open caps at infinity $\BH^2\times\{\pm\infty\}$ in \cite{KM} to study infinite curves. To state their result, and clarify the ambiguity, we introduce the following notation.

\begin{defn} \label{cornerdefn} Let $\Gamma_g^\pm=\Gamma\cap \BH^2\times\{\pm\infty\}$ corresponds to the interior of $\Gamma^\pm$. In particular, if nonempty, $\Gamma_g^\pm$ will be a collection of arcs in $\BH^2\times\{\pm\infty\}$. The subscript $g$ corresponds to the term "geodesic", which will be clear in the next lemma.

Let $\overline{\Gamma_g^\pm}$ be the closure of $\Gamma_g^\pm$ in $\overline{\BH^2}\times\{\pm\infty\}$. In other words, if $\Gamma_g^\pm$ is a finite collection of arcs in $\BH^2\times\{\pm\infty\}$, then $\overline{\Gamma_g^\pm}$ is $\Gamma_g^\pm$ with the endpoints in $S^1_\infty\times\{\pm\infty\}$.

Let $\Gamma_c^\pm=\Gamma^\pm-\overline{\Gamma_g^\pm}$ be the remaining boundary points of $\Gamma^\pm$ in $S^1_\infty\times\{\pm\infty\}$. Here, the subscript $c$ correspond o the term "corner". In particular, $\Gamma_c^\pm$ are exactly the points in $\Gamma$ which are not transverse to the corner circles $S^1_\infty\times\{\pm\infty\}$. In other words, $\Gamma_c^\pm$ are the points in $S^1_\infty\times\{\pm\infty\}$ which are interior points of the arcs of $\Gamma\cap (S^1_\infty\times \overline{\BR})$.

We will call $\Gamma$ {\em nonoverlapping at the corner} if $\Gamma_c^\pm$ does not contain any interval in the corner circles $\CS$. In other words, $\Gamma$ does not overlap with the corner circles of $\SI$ in any interval. This property will be used in the next section to classify strongly fillable curves.
\end{defn}

Now, we can state the key property of infinite fillable curves, given in \cite[Proposition 4.3]{KM}.

\begin{lem} \label{geod} \cite{KM} If $\Gamma$ is an infinite fillable curve in $\SI$, then $\Gamma_g^\pm$ must be a collection of geodesics in $\caps$.
\end{lem}


This lemma is crucial to understand the structure of $\Gamma^\pm$ for fillable curves. The following lemma will also be useful in the following sections.

\begin{lem} \label{geodlem2} If $\Gamma$ is an infinite fillable curve in $\SI$, then $\overline{\Gamma_g^\pm}$ in $\overline{\BH^2}\times\{\pm\infty\}$ must be disjoint.
\end{lem}

\begin{pf} By Lemma \ref{geod}, we only need to show that two geodesics in $\Gamma_g^+$ cannot have same endpoint in $S^1_\infty\times\{+\infty\}$. Assume on the contrary that $\gamma_1$ and $\gamma_2$ in $\Gamma^+$ have same endpoint $p=(\theta_0,\infty) \in S^1_\infty\times\{+\infty\}$. Let $\R=[\theta_1,\theta_2]\times [c_1,c_2]$ be a rectangle in $\Si$ such that $\theta_0\in (\theta_1,\theta_2)$ and $\R\cap \Gamma=\emptyset$. Let $\T$ be the unique area minimizing surface in $\BHH$ with $\PI \T=\partial \R$. Let $\Omega$ be the component of $\BHH-\T$ with $\PI\Omega=int(\R)$. By \cite[Lemma 2.11]{Co}, we can foliate $\Omega$ by unique area minimizing surfaces bounding smaller rectangles in $\R$.

By the proof of \cite[Proposition 4.3]{KM}, the limit of $\Sigma_t=\Sigma-t$ is a minimal plane $\wh{\Sigma}$ which is a product plane, $\Gamma^+_g\times\BR$. By construction, $\Gamma_t\cap\R=\emptyset$ for any $t$. As $\Omega$ is foliated by minimal planes, this implies $\Sigma_t\cap \Omega=\emptyset$ by maximum principle. However, the limit plane $\wh{\Sigma}\cap\Omega\neq \emptyset$ as $\wh{\Sigma}\supset (\gamma_1\cup\gamma_2)\times \BR$. This is a contradiction. The proof follows.
\end{pf}

Now, we give a natural generalization of tall rectangles to infinite curves.

\begin{defn} \label{infrecdef} \cite{ST} [Infinite Rectangles] Let $\gamma$ be a complete geodesic in $\BH^2$ with $\PI\gamma=\{p,q\}$. Let $\alpha$ be one of the two arcs in $S^1_\infty(\BH^2)$ with endpoints $p$ and $q$. Fix $t_0\in\BR$. Let $l^+_p=\{p\}\times[t_0,\infty]$ and $l^-_p=\{p\}\times[-\infty,t_0]$ be the vertical line segments in $\Si$. Let $\gamma^\pm=\gamma\times\{\pm\infty\}$ be the geodesic in $\caps$. Let $\alpha_0=\alpha \times\{t_0\}$. Then define $\R^+=\gamma^+\cup l^+_p\cup l^+_q\cup \alpha_0$ is an infinite rectangle. Similarly, $\R^-=\gamma^-\cup l^-_p\cup l^-_q\cup \alpha_0$ is also an infinite rectangle. See Figure \ref{fig2}-Left, where each component is an infinite rectangle.
\end{defn}


\begin{lem} \label{infrec} Any infinite rectangle in $\SI$ is strongly fillable, and it bounds a unique area minimizing surface in $\SI$.
\end{lem}

\begin{pf} Let $\R=\gamma^+\cup l^+_p\cup l^+_q\cup\alpha_0$ be an infinite rectangle with the notation above. By \cite{ST}, $\R$  bounds a minimal surface $\T$ in $\BHH$, i.e. $\PI \T=\R$. Furthermore, $\T$ is a graph over the region $\Delta$ in $\overline{\BH^2}$ separated by $\gamma\cup \alpha$ \cite{ST}.

Consider the family of minimal surfaces $\{\T_s\mid s\in \BR\}$ where $\T_s$ is $s$ vertical translation of $\T$. By construction, $\{\T_s\}$ foliates the convex region $\Delta\times \BR$. This shows that $\T$ is area minimizing, and the unique minimal surface $\R$ bounds in $\BHH$. Hence, $\R$ is strongly fillable.
\end{pf}

\begin{defn}\label{tamedefn} [Tame Curves]
We will call an infinite curve $\Gamma$ in $\SI$ {\em tame} if $\Gamma^\pm$ has finitely many component. Otherwise, we will call $\Gamma$ {\em a wild curve.} 
\end{defn}

\begin{rmk} Throughout the paper, all the curves are assumed to be tame unless otherwise stated. In Section \ref{wild}, we give two Cantor-like examples of wild curves which point out that the asymptotic Plateau problem can be quite different for wild curves in general. 
\end{rmk}

\subsection{Scherk Graphs} \label{scherksec} \

\begin{figure}[h]
\begin{center}
$\begin{array}{c@{\hspace{.4in}}c}

\relabelbox  {\epsfysize=2in \epsfbox{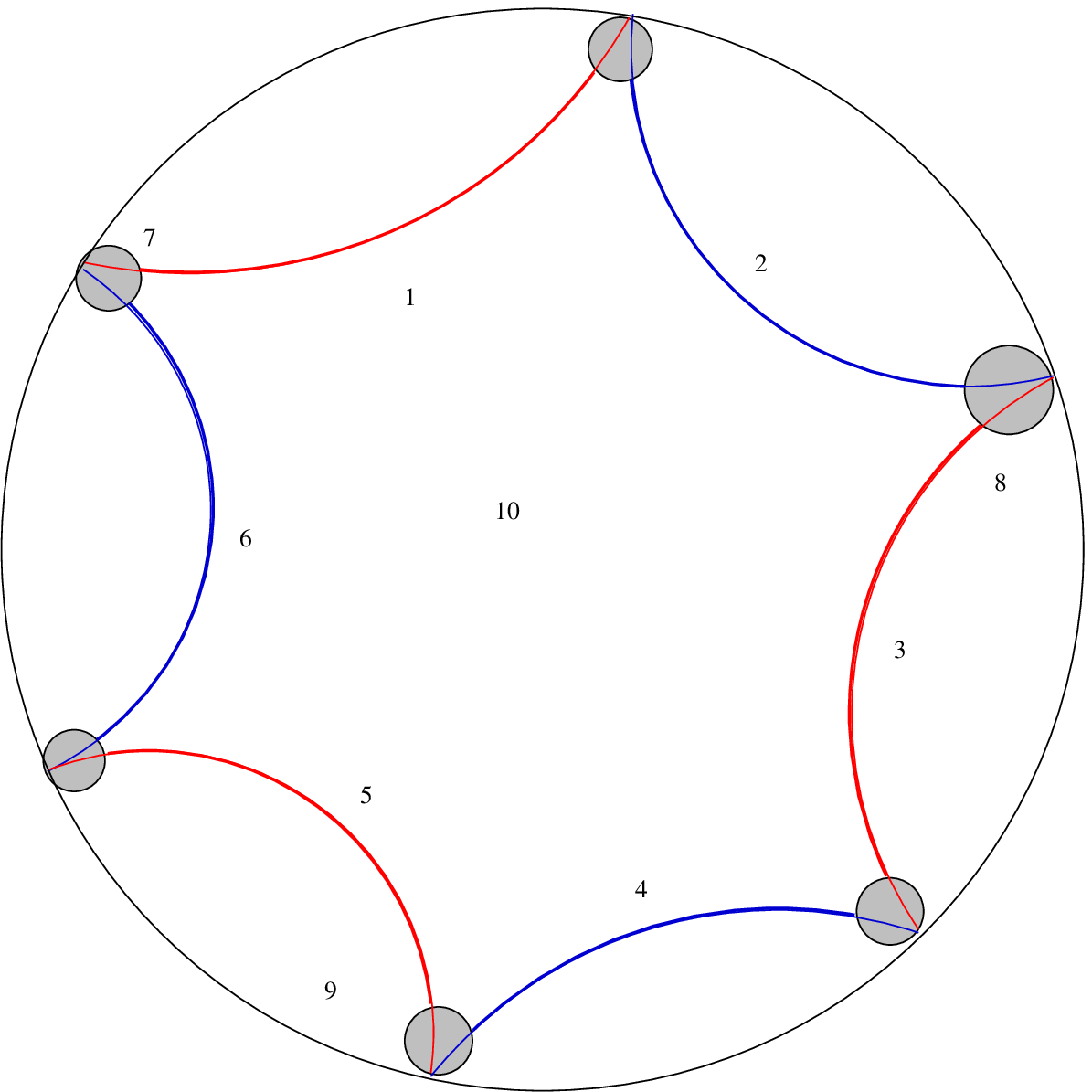}}
\relabel{1}{\small $\wh{\alpha}_1$} \relabel{3}{\small $\wh{\alpha}_3$} \relabel{5}{\small $\wh{\alpha}_2$}

\relabel{2}{\small $\wh{\beta}_3$} \relabel{4}{\small $\wh{\beta}_2$} \relabel{6}{\small $\wh{\beta}_1$}

\relabel{7}{\footnotesize $B_2$} \relabel{8}{\footnotesize  $B_4$} \relabel{9}{\footnotesize $B_6$} \relabel{10}{$\Delta$} \endrelabelbox &

\relabelbox  {\epsfysize=2in \epsfbox{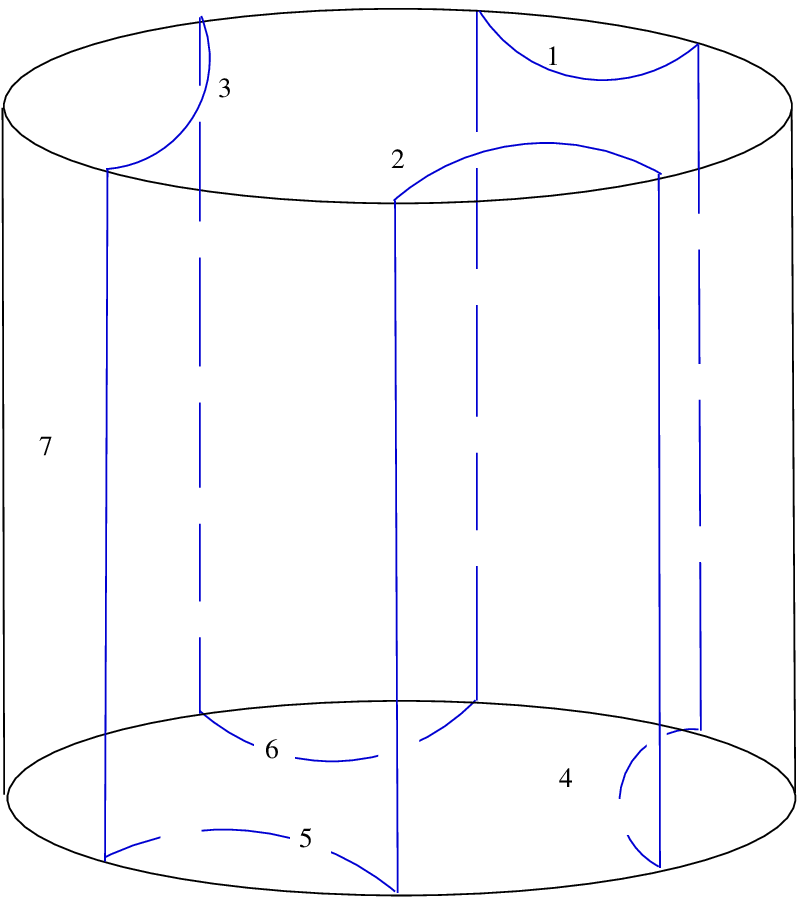}}
\relabel{1}{\scriptsize  $\alpha^+_3$} \relabel{2}{\scriptsize  $\alpha^+_2$ }  \relabel{3}{\scriptsize $\alpha^+_1$}
\relabel{4}{\scriptsize $\beta^-_2$} \relabel{5}{\scriptsize $\beta^-_1$ }  \relabel{6}{\scriptsize $\beta^-_3$}
\relabel{7}{\footnotesize $\Gamma$}  \endrelabelbox \\
\end{array}$

\end{center}

\caption{ \label{ideal-scherk} \footnotesize In the figure left, $\Delta$ is an ideal hexagon where $B_i$ represents the small horoballs at $p_i\in S^1_\infty$. In the right, the asymptotic boundary of a Scherk graph is given.}
\end{figure}

Now, we recall the results on Scherk graphs in $\BHH$ by \cite{CR}. These are minimal graphs over ideal $2n$-gons in $\BH^2$ where the graph takes values $+\infty$ and $-\infty$ on alternating sides. In particular, let $\Delta$ be a closed ideal $2n$-gon in $\BH^2$. Let $\mathcal{V}=\{p_1,p_2, ..., p_{2n}\}\in \si$ be the set of ideal vertices of $\Delta$ which are circularly ordered. Let $\alpha_i$ be the geodesic with $\PI\alpha_i = \{p_{2i-1},p_{2i}\}$ and $\beta_i$ be the geodesic with $\PI\beta_i = \{p_{2i},p_{2i+1}\}$. Then, $\partial \Delta = \alpha_1\cup\beta_1\cup...\cup\alpha_n\cup\beta_n$. For each ideal vertex $p_i$, define a sufficiently small horocycle $C_i$ such that $C_i\cap C_j=\emptyset$ for any $1\leq i<j\leq 2n$. Let  $B_i$ be the open horodisk which $C_i$ bounds in $\BH^2$. Let $\Omega=\BH^2-\bigcup_i B_i$. Let $\wh{\alpha}_i=\alpha_i\cap\Omega$ and $\wh{\beta}_i=\beta_i\cap\Omega$. Let $a(\Delta)=\Sigma |\wh{\alpha}_i|$ and $b(\Delta)=\Sigma |\wh{\beta}_i|$. See Figure \ref{ideal-scherk}-left.

We say an ideal polygon $\D$ is {\em inscribed in} $\Delta$ if $\PI\D\subset\mathcal{V}$. Clearly, $\partial \D$ consists of some geodesics in $\partial \Delta$, and some other geodesics $\{\gamma_j\}$ in the interior of $\Delta$. Let $a(\D)$ is the sum of $|\wh{\alpha}_i|$ where $\alpha_i\subset \partial \D$, and similarly define $b(\D)$. Let $c(\D)=\Sigma |\wh{\gamma_j}|$ where $\wh{\gamma}_j=\gamma_j\cap\Omega$. Then, let $|\D|$ be the sum of the "truncated lengths" of the geodesics in $\partial \D$, i.e. $|\D|=a(\D)+b(\D)+c(\D)$.

\begin{defn} [Exact Polygons] \label{exact} Let $\Delta$ be an ideal $2n$-gon in $\BH^2$. For any inscribed polygon $\D$ in $\Delta$ with $\D\neq \Delta$, let $2a(\D)<|\D|$ and $2b(\D)<|\D|$. Then, we call $\Delta$ a {\em regular ideal polygon}.

Let $\Delta$ be a regular ideal polygon with $a(\Delta)=b(\Delta)$. Then, we call $\Delta$ an {\em exact ideal polygon}.
\end{defn}

Recently, Collin and Rosenberg showed the existence of solutions to the Dirichlet problem with $\pm\infty$ boundary values for exact ideal polygons \cite[Theorem 1]{CR}.

\begin{lem} \label{scherk} \cite{CR} [Scherk Graphs] Let $\Delta$ be an exact ideal $2n$-gon in $\BH^2$. Then, there exists a solution $u:\Delta\to \BR$ to the minimal surface equation on $\Delta$ which takes values $+\infty$ on $\alpha_i$ and $-\infty$ on $\beta_i$ for $1\leq i\leq n$. Furthermore, the solution is unique up to an additive constant.
\end{lem}

\begin{rmk} \label{scherkrem} [Area Minimizing] The Scherk graph $\Sigma=graph(u)$ is a minimal surface in $\BH^2$ with $$\xi=\PI \Sigma = \bigcup_{i=1}^n (\alpha_i\times\{+\infty\}) \cup (\beta_i\times\{-\infty\})\bigcup_{j=1}^{2n} l_{p_j}$$ where $l_p$ is the vertical line $\{p\}\times\BR$ in $\Si$. We will call the asymptotic boundary $\xi$ of a Scherk graph a {\em Scherk curve} in $\SI$. See Figure \ref{ideal-scherk}-right.

Notice that $\Sigma$ is also an area minimizing surface in $\BHH$ because the family of surfaces $\{\Sigma_t\mid t\in \BR\}$ foliates the convex region $\Delta\times \BR$ where $\Sigma_t$ is the $t$ vertical translation of $\Sigma$.
\end{rmk}

\subsection{Fat / Skinny at Infinity} \label{fatsec} \

Now, we introduce a new notion to study the fillability of the infinite curves by using Scherk graphs. Let $\Gamma=\Gamma^+\cup\Gamma^-\cup\wt{\Gamma}$ as before. Since we assume $\Gamma$ is tame, $\Gamma^\pm_g$ is a finite collection of geodesics. Let $\Gamma^+_g=\gamma_1^+\cup ...\gamma^+_n$ and $\Gamma^-_g=\gamma_1^-\cup ...\gamma^-_m$ for $n,m>1$ where $\gamma_i^\pm$ corresponds to a geodesic in $\BH^2\times\{\pm\infty\}$ by Lemma \ref{geod}. $n=1$ or $m=1$ cases are trivial, and they will be discussed later. 

By Lemma \ref{geodlem2}, we further assume no endpoints of $\gamma^+_i$ and $\gamma^+_j$ are the same for $i\neq j$. Let $\V^+=\bigcup \PI\gamma^+_i$ be the set of $2n$ points in $S^1_\infty\times\{+\infty\}$. Let $\V^+= \{p_1^+,p_2^+,...,p_{2n}^+\}$ be indexed so that the points are circularly ordered. For $1\leq i\leq 2n$, let $\tau^+_i$ be the geodesic in $\BH^2\times\{+\infty\}$ with $\PI \tau_i^+= \{p_i^+,p_{i+1}^+\}$. Of course, for some $i$, $\tau^+_i\subset \Gamma^+$. Similarly, define $\V^-$ to be the set of endpoints of $\gamma^-_j$, and $\tau^-_j$ to be the geodesics between them for $1\leq j\leq 2m$.

\begin{figure}[h]
	\begin{center}
		$\begin{array}{c@{\hspace{.4in}}c}

		\relabelbox  {\epsfysize=2in \epsfbox{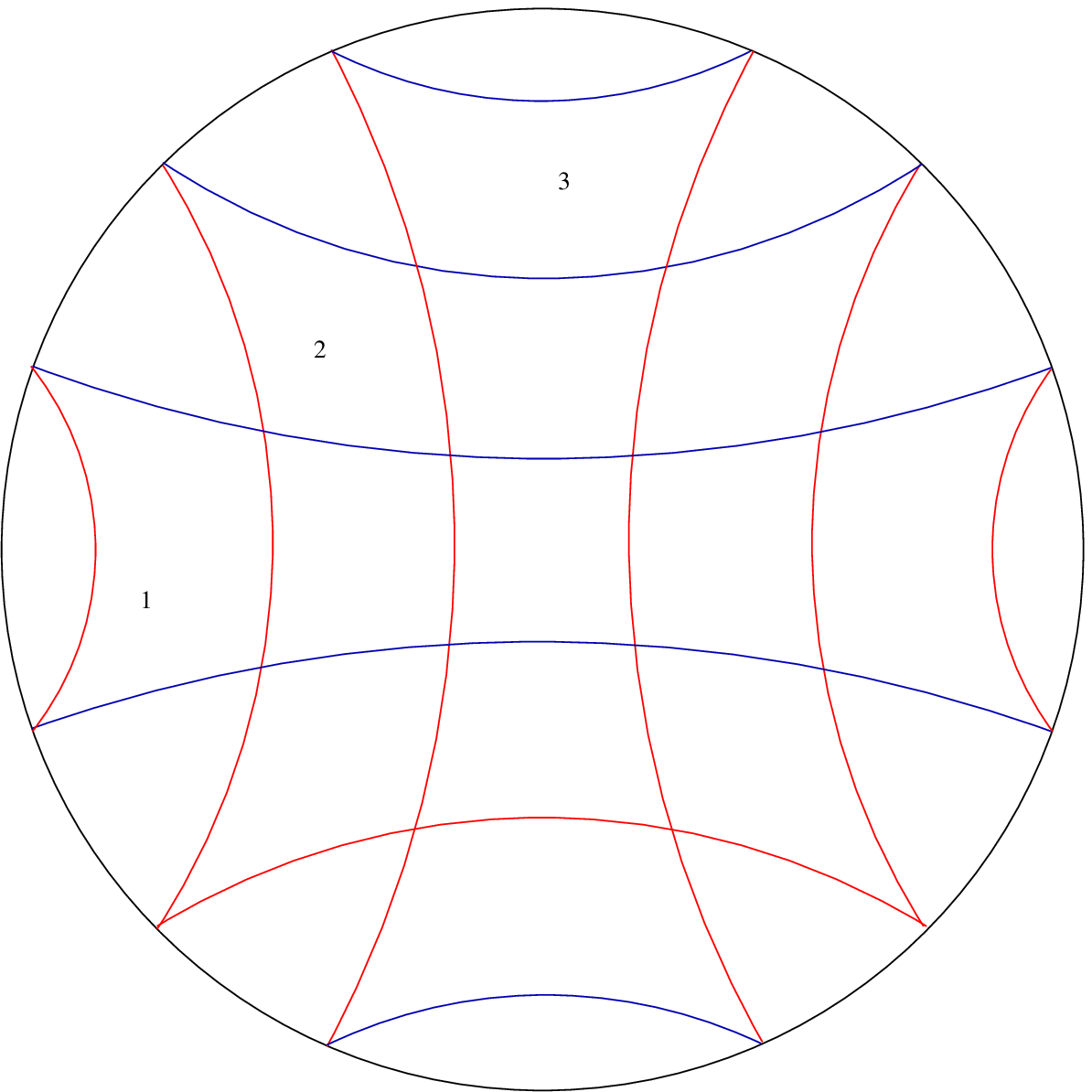}} \relabel{1}{$\Delta_1$} \relabel{2}{$\Delta_2$} \relabel{3}{$\Delta_3$} \endrelabelbox &
		
		\relabelbox  {\epsfysize=2in \epsfbox{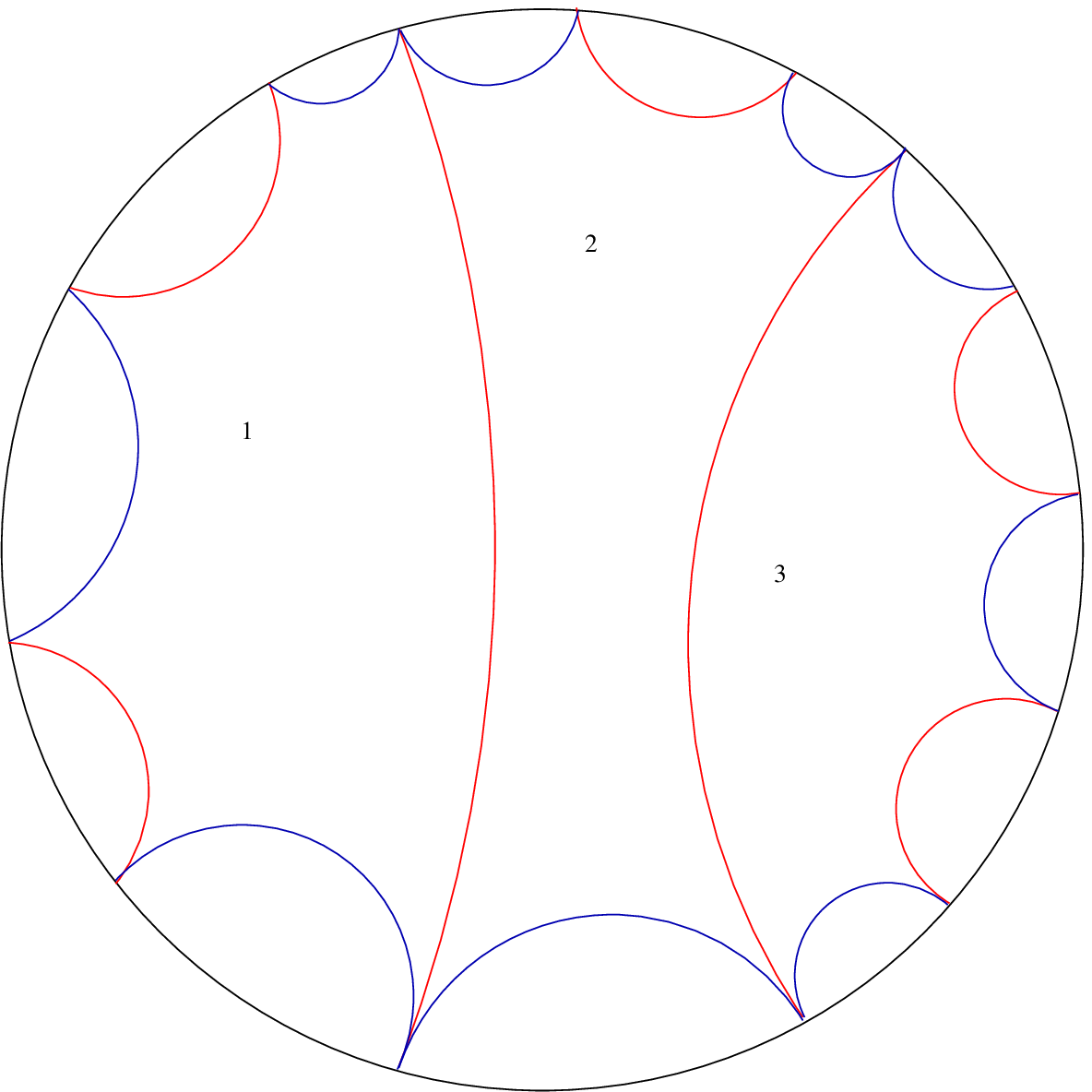}}
		\relabel{1}{$\Delta^+_1$} \relabel{2}{$\Delta^+_2$ }  \relabel{3}{$\Delta^+_3$}   \endrelabelbox \\
		\end{array}$
		
	\end{center}
	
	\caption{ \label{polygon_fig} \footnotesize In the figure left, $\Delta_1,\Delta_2$ and $\Delta_3$ represent a fat, exact, and skinny polygon respectively. Red curves represent $\alpha$ curves, and blue curves represent $\beta$ curves in $\partial\Delta_i$. In the figure right, $7$ red geodesics represent $\Gamma^+$, and $\Delta^+$ decomposed into $3$ inscribed polygons   $\Delta^+ = \Delta^+_1 \cup \Delta^+_2 \cup \Delta^+_3$.}
\end{figure}

Let $\Delta^+$ be the convex hull of $\V^+$ in $\BH^2\times\{+\infty\}$. In particular, $\Delta^+$ is the closed ideal $2n$-gon in $\BH^2\times\{+\infty\}$ with $\partial \Delta^+=\bigcup_{i=1}^{2n} \tau^+_i$. Similarly define $\Delta^-$ which is a closed ideal $2m$-gon with $\partial \Delta^-=\bigcup_{j=1}^{2m} \tau^-_j$. Notice  that some of the geodesics $\gamma^+_i$ might be in the interior of $\Delta^+$. Hence, $\Delta^+-\Gamma^+$ is a union of ideal polygons in $\BH^2\times\{+\infty\}$, i.e. $\Delta^+-\Gamma^+=\Delta^+_1\cup..\cup\Delta^+_{n_1}$ for some $n_1\geq 1$, where the vertices of $\Delta_i^+$ is in $\V$ (See Figure \ref{polygon_fig}-right). Similarly, let $\Delta^- -\Gamma^-=\Delta^-_1\cup...\Delta^-_{m_1}$. By abuse of notation, we will take $\Delta^\pm_i$ as the closed polygons containing its sides, i.e. $\Delta^+=\Delta^+_1\cup...\cup\Delta^+_{n_1}$ and $\Delta^- =\Delta^-_1\cup..\cup\Delta^-_{m_1}$. Notice that each $\Delta^\pm_i$ is inscribed in $\Delta^\pm$.

In particular, $\Delta^\pm$ naturally decomposes as a union of {\em inscribed} ideal polygons $\Delta^\pm_i$. Notice that each $\partial\Delta^\pm_i$ contains $k_i^\pm$ geodesics from $\Gamma^\pm$ by construction, and hence, $\Delta^\pm_i$ is $2k^\pm_i$-gon for some $k^\pm_i>1$. Furthermore, the geodesics in $\Gamma^+$ alternates in $\partial\Delta^+_i$, and similarly $\Gamma^-$ in $\partial\Delta^-_j$. Declare $\Gamma^+$ in $\partial \Delta^+_i$ as $\alpha$-curves of $\Delta^+_i$. Similarly for $\Delta^-_j$. Hence, if $\Delta^\pm$ is regular, then so is $\Delta^\pm_i$  as $\Delta^\pm_i$ is inscribed polygon in $\Delta^\pm$ with the induced $\alpha$ and $\beta$ curves.

\begin{defn} [Fat / Skinny at Infinity] \label{fatdef} For $n>1$, let $\Omega$ be a regular ideal $2n$-gon in $\BH^2$ with $\partial\Omega=\alpha_1\cup\beta_1\cup...\alpha_n\cup\beta_n$ where the geodesics $\alpha_i$  and $\beta_i$ are alternating. Let $a(.)$ and $b(.)$ be as defined in Scherk graphs section. We will call $\Omega$ {\em fat} if $a(\Omega)<b(\Omega)$. We will call $\Omega$ {\em skinny} if $a(\Omega)>b(\Omega)$. See Figure \ref{polygon_fig}-left.

Let $\Gamma$ be a tame infinite curve in $\SI$. Define $\Delta^\pm$ as above. Let $\Delta^+=\Delta^+_1\cup...\cup\Delta^+_{n_1}$ and $\Delta^- =\Delta^-_1\cup..\cup\Delta^-_{m_1}$ be induced decompositions of $\Delta^\pm$ as defined before. Let the geodesics in $\Gamma^+\cap\partial \Delta_i^+$ be $\alpha$-curves of $\Delta_i^+$, and let the geodesics in $\Gamma^-\cap\partial \Delta_j^-$ be again $\alpha$-curves of $\Delta_j^-$. Then, we will call $\Gamma$ {\em fat at infinity} if all ideal polygons $\Delta^+_i$ and $\Delta^-_j$ are fat. Furthermore, we will call $\Gamma$ {\em skinny at infinity} if at least one polygon $\Delta^\pm_i$ is skinny. 
\end{defn}

Note that when we say $\Gamma$ fat at infinity, we implicitly assume that $int(\Gamma^\pm)$ is already a collection of geodesics in $\caps$. In the case $int(\Gamma^\pm)=\emptyset$, we also say $\Gamma$ is fat at infinity.Note that if $n=1$ or $m=1$ (number of components in $\Gamma^\pm$), we still call $\Gamma$ is fat at infinity, and the proofs in Section \ref{sfillsec} apply to this case trivially. 

\begin{rmk} [Symmetric definition in $\Delta^-$] Notice that we are using a symmetry between $\Gamma^+$ and $\Gamma^-$ in the definition above. Alternatively, one can define {\em being fat at infinity} as follows. We declare $\Gamma^+_g$ curves as $\alpha$ curves in $\Delta_i^+$ and require $a(\Delta^+_i)<b(\Delta^+_i)$ to call $\Delta_i^+$ fat as before. On the other hand, we can call $\Gamma^-_g$ curves as $\beta$ curves in $\Delta_i^-$ and require $b(\Delta^+_i)<a(\Delta^+_i)$ to call $\Delta_i^-$ fat. Because of the Scherk graphs we are going to use, this perspective might seem more natural. In any case, it is not hard to see that both definitions are equivalent.
\end{rmk}

The following lemma implies that any fat polygon can be covered by a finite union of exact polygons. See Figure \ref{fatcovering}. In particular, let $\Omega$ be a fat $2n$-gon with the vertices $\V=\{p_1, ...p_{2n}\}$ which is circularly ordered. Hence, $\Omega$ is the convex hull of $\V$ in $\BH^2$. Define $p_{2n+1}=p_1$.  Let $\alpha_i=\overline{p_{2i-1}p}_{2i}$ and $\beta_i=\overline{p_{2i}p}_{2i+1}$ where $\overline{pq}$ represents the geodesic between $p$ and $q$. Hence, $\partial \Omega = \alpha_1\cup\beta_1\cup..\cup\alpha_n\cup\beta_n$ where the geodesics $\alpha_i$ and $\beta_i$ are alternating.

\begin{lem} \label{fatlem} [Fats covered by Exacts] Let $\Omega$ be a fat $2n$-gon as above. Then, $\Omega$ can be covered by exact $2n$-gons $\D_i$ with $\Omega\subset\bigcup_i\D_i$ such that $int(\D_i)\cap\alpha_j=\emptyset$ for any $i,j$. 
\end{lem}

\begin{pf}  For $p,q\in \si$, let $[p,q]$  represent the interval from $p$ to $q$ in $\si$ in the counterclockwise direction. We represent counterclockwise circular order of $\V$ with $p_1\prec p_2 \prec ...\prec p_{2n}\prec p_1$. 
	

\begin{figure}[h]
\begin{center}
$\begin{array}{c@{\hspace{.4in}}c}

\relabelbox  {\epsfysize=2in \epsfbox{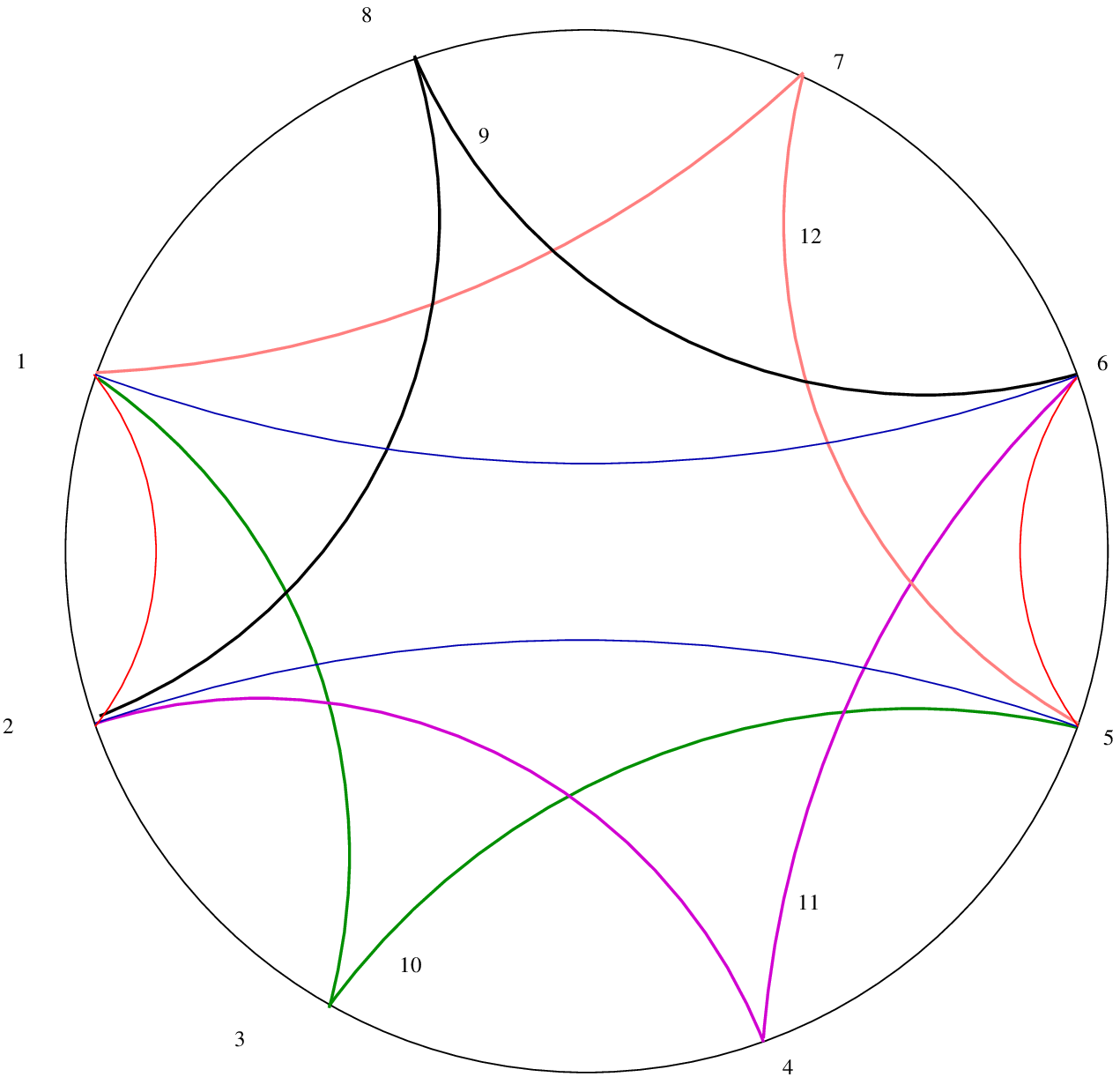}}
\relabel{1}{\footnotesize $p_1$} \relabel{2}{\footnotesize $p_2$ }  \relabel{3}{\footnotesize $p_2^*$} \relabel{4}{\footnotesize $p_3^*$} \relabel{5}{\footnotesize $p_3$} \relabel{6}{\footnotesize $p_4$ } 
\relabel{7}{\footnotesize $p_4^*$} \relabel{8}{\footnotesize $p_1^*$}
\relabel{9}{\footnotesize $\D_1$} \relabel{10}{\footnotesize $\D_2$ }  \relabel{11}{\footnotesize $\D_3$} \relabel{12}{\footnotesize $\D_4$}
\endrelabelbox  &

\relabelbox  {\epsfysize=2in \epsfbox{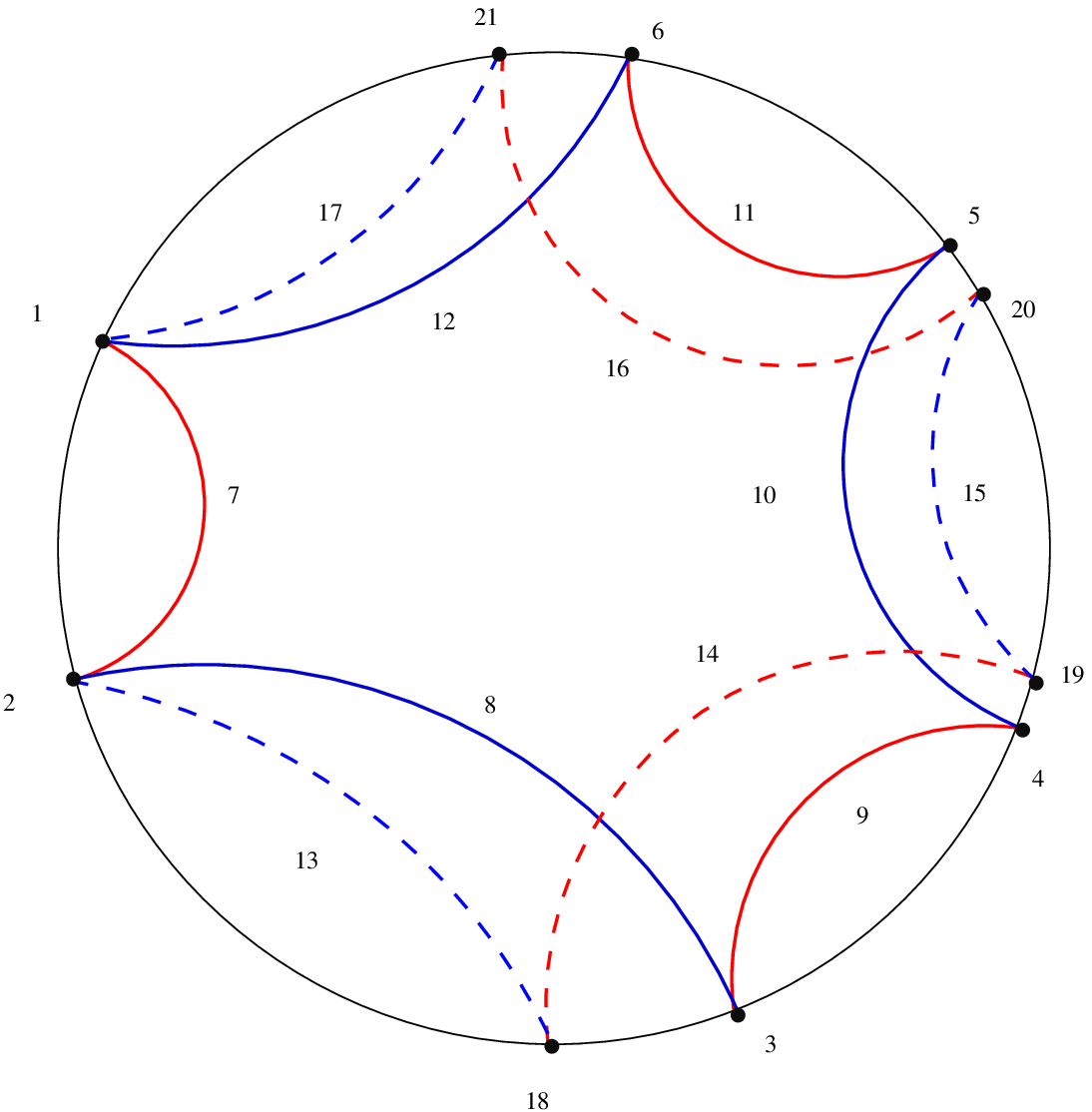}}
\relabel{1}{\footnotesize $p_1$} \relabel{2}{\footnotesize $p_2$ }  \relabel{3}{\footnotesize $p_3$} \relabel{4}{\footnotesize $p_4$} \relabel{5}{\footnotesize $p_5$} \relabel{6}{\footnotesize $p_6$ } 

\relabel{7}{\footnotesize $\alpha_1$} \relabel{8}{\footnotesize $\beta_1$} \relabel{9}{\footnotesize $\alpha_2$} \relabel{10}{\footnotesize $\beta_2$ }  \relabel{11}{\footnotesize $\alpha_3$} \relabel{12}{\footnotesize $\beta_3$}

\relabel{13}{\footnotesize $\beta_1'$} \relabel{14}{\footnotesize $\alpha_2'$} \relabel{15}{\footnotesize $\beta_2'$ }  \relabel{16}{\footnotesize $\alpha_3'$} \relabel{17}{\footnotesize $\beta_3'$}

\relabel{18}{\footnotesize $p_3^*$} \relabel{19}{\footnotesize $p_4^*$} \relabel{20}{\footnotesize $p_5^*$} \relabel{21}{\footnotesize $p_6^*$ }

\endrelabelbox \\
\end{array}$

\end{center}

\caption{ \label{fatcovering} \footnotesize [Exact Covering] In the figure left, $\Omega$ represents a fat $4$-gon with vertices $p_1,p_2,p_3,p_4$. The exact covering is the one given in the Lemma \ref{fatlem}. In the figure right, $\Omega$ is a fat hexagon, and $\D_1$ is exact hexagon described with dashed lines (Special covering in Remark \ref{specialcover}).}
\end{figure}

Let $i\in\{1,2,...,2n\}$ and choose $p^*_i$ in $(p_{i-1},p_i)$ if $i$ is odd, and 
choose $p^*_i$ in $(p_i,p_{i+1})$ if $i$ is even. Let $\W_i=\{p_1, p_2, .., p_{i-1}, p^*_i, p_{i+1},..., p_{2n}\}$ such that $\V\triangle\W_i=\{p_i,p^*_i\}$. Let $\D_i$ be the convex hull of $\W_i$. We claim that we can choose a unique such $p_i^*$ so that $\D_i$ is exact (See Fig \ref{fatcovering}-left).


We describe the exact covering as follows: For $i=2k$, $\D_i$ has the same geodesics with $\Omega$ except $\alpha_k$ and $\beta_k$. Intuitively, $\alpha_k$ gets "bigger" while $\beta_k$ gets "smaller". In particular, for $i=2k$, if $p^*_i\to p_{i+1}$, then $\beta_k$ escapes to infinity, and the quantity $b(\D_i)-a(\D_i)\searrow -\infty$ monotonically as $p^*_i\to p_{i+1}$. Also, if $p^*_i$ moves to other direction, $p^*_i\to p_{i}$, then $\D_i\to\Omega$ and $b(\D_i)-a(\D_i)\nearrow (b(\Omega)-a(\Omega))>0$ monotonically. This proves that there exists a unique point $p^*_i\in(p_i,p_{i+1})$ with $a(\D_i)=b(\D_i)$.

For $i=2k+1$, $\D_i$ has the same geodesics with $\Omega$ except $\alpha_{k+1}$ and $\beta_k$. This time intuitively, $\alpha_{k+1}$ gets "bigger" while $\beta_k$ gets "smaller". Similar to the even case, if $i=2k+1$, when $p^*_i\to p_{i-1}$, $\beta_k$ escapes to infinity, and the quantity $b(\D_i)-a(\D_i)\searrow -\infty$ monotonically. In the other direction, if $p^*_i\to p_{i}$, then $\D_i\to\Omega$ and $b(\D_i)-a(\D_i)\nearrow (b(\Omega)-a(\Omega))>0$ monotonically again. Hence, there exists a unique point $p^*_i\in (p_{i-1},p_i)$ with $a(\D_i)=b(\D_i)$ in this case, too.

Notice that for both cases, the new geodesics $\alpha_k^*$ and $\beta_k^*$ intersect only $\beta_k$, and they don't intersect $\alpha_k$. Hence, $int(\D_i)\cap \alpha_j=\emptyset$ for any $j$.

Assuming $\D_i$ is a regular polygon for $1\leq i\leq 2n$, we finish the proof as follows. Each $\D_i$ is an exact $2n$-gon which has the same sides with $\Omega$ except $\beta_k$ and $\alpha_k$  for $i=2k$ ($\alpha_{k+1}$ for $i=2k+1$). Then, it is not hard to show that $\Omega\subset \bigcup_{i=1}^{2n}\D_i$. Let $\V_i=\V-\{p_i\}$ for $1\leq  i\leq 2n$. $\wt{\D}_i$ be the ideal ($2n-1$)-gon with vertices $\V_i$. As $\V_i\subset \W_i$, then $\wt{\D}_i\subset \D_i$. In particular, $\wt{\D}_i$ is obtained by removing the ideal triangle $\Delta_i$ with vertices $\{p_{i-1},p_i,p_{i+1}\}$ from $\Omega$, i.e. $\wt{\D}_i=\Omega-\Delta_i$. Since $\wt{\D}_i\subset \D_i$ for any $1\leq i\leq 2n$, it is clear that $\Omega\subset\bigcup_{i=1}^{2n}\D_i$, and the proof follows with the following claim.

\vspace{.2cm}

\noindent {\bf Claim:} $\D_i$ is a regular polygon.

\vspace{.2cm}

\noindent {\em Proof of the Claim:} Let $E$ be an inscribed polygon in $\D_i$. Let $v(E)$ be the vertices of $E$. Hence, $v(E)\subset \W_i$, the vertices of $\D_i$. We need to show that $a(E)<b(E)+c(E)$ and $b(E)<a(E)+c(E)$ with the notation in the section \ref{scherksec}.

First assume that $p_i^*$ is not in $v(E)$. Then, $E$ is also an inscribed polygon in $\Omega$ which is regular. So, $E$ satisfies the inequalities.

Now, assume that $p_i^*\in v(E)$. There are two cases. The first case is that $p_i^*$ is not an endpoint of a $\gamma$ curve of $E$. Let $E'$ be the ideal polygon with the same vertices with $E$ except $p_i$ replaced with $p_i^*$. Since $E'$ is inscribed polygon in $\Omega$, it satisfies the inequalities. Notice that by construction, $|\wh{\alpha}_k|-|\wh{\beta}_k|<|\wh{\alpha^*}_k|-|\wh{\beta^*}_k|$. As the other sides coincides with $E'$, $b(E)<a(E)+c(E)$ follows. Now, we will show $a(E)<b(E)+c(E)$. Notice that $E^c$, the complement of $E$ in $\D_i$, is an inscribed polygon in $\Omega$ as $p_i^*$ is not in $v(E^c)$. Hence, $E^c$ satisfies the inequalities, i.e. $b(E^c)<a(E^c)+c(E^c)$. As $\D_i$ is exact, $a(\D_i)=b(\D_i)=\mu$. Notice that $b(\D_i)=b(E)+b(E^c)$ and $a(\D_i)=a(E)+a(E^c)$. Also, $c(E)=c(E^c)$. Hence, $b(E^c)<a(E^c)+c(E^c)$ implies that $(\mu-b(E))<(\mu-a(E))+c(E)$. Hence, we obtain $a(E)<b(E)+c(E)$

Now, for the second case, assume that $p_i^*$ is an endpoint of a $\gamma$ curve $\gamma^*_o$ of $E$. Assume $\beta_k^*$ is in $E$. Let the other endpoint of $\gamma^*_o$ be $p_{j_o}$. Let $\gamma_o$ be the geodesic with the endpoints $p_{j_o}$ and $p_i$. Let $E'$ be the ideal polygon with the same vertices with $E$ except $p_i$ replaced with $p_i^*$. Since $E'$ is inscribed polygon in $\Omega$, it satisfies the inequalities. $E$ and $E'$ has the same sides except that $E'$ has $\beta_k$ and $\gamma_o$ while $E$ has the sides $\beta_k^*$ and $\gamma_o^*$. Consider the hyperbolic isometry $\varphi$ fixing $p_{j_o}$ and $p_{i+1}$, and sending $p_i$ to $p_i^*$. Then, $\varphi(\gamma_o)=\gamma_o^*$ and $\varphi(\beta_k)=\beta_k^*$. Hence, $|\wh{\gamma}_0|-|\wh{\beta}_k|=|\wh{\gamma^*}_o|-|\wh{\beta^*}_k|$. As $E'$ satisfies the inequalities, this implies $E$ satisfies the inequalities, too. Similarly, if $E$ contains $\alpha_k^*$ instead of $\beta_k^*$ same argument would work by replacing the hyperbolic isometry $\varphi$  with $\psi$ fixing $p_{j_o}$ and $p_{i-1}$ (the other endpoint of $\alpha_k$), and sending $p_i$ to $p_i^*$. The proof follows.
\end{pf}


We will call such coverings of fat $2n$-gon $\Omega$ by exact $2n$-gons $\{\D_i\}$ as in the lemma, an {\em exact covering of} $\Omega$. See Figure \ref{fatcovering}-left.

\begin{rmk}  [Special Exact Covering] \label{specialcover} In previous lemma, we give a simple way to cover a fat polygon by exact polygons. In the following, we will need a modified version of this covering as follows. For a given fat $2n$-gon $\Omega$ with $\partial \Omega = \alpha_1\cup\beta_1\cup..\cup\alpha_n\cup\beta_n$, we will give a covering with $n$ exact polygons $\D_n$. For any fixed $1\leq i_o\leq n$, let $\D_{i_o}$ be the exact polygon containing $\alpha_{i_o}$, and for any $k\neq i_o$, $\alpha_k$  gets "bigger",  while $\beta_i$ gets "smaller" for any $1\leq i\leq n$. In particular, we fix the endpoints of $\alpha_{i_o}$, i.e. $p_{2i_o-1}$ and $p_{2i_o}$. However, we move all other odd indexed $p_{2k-1}$ clockwise a little bit, while we move all other even indexed $p_{2k}$ counterclockwise a little bit (See Figure \ref{fatcovering}-right). As $\Omega$ is fat polygon, when we move the remaining points enough amount, this will give us an exact polygon $\D_{i_o}$ containing $\alpha_{i_o}$. Then, it is not hard to see we will get a similar exact covering of $\Omega$ by $\{\D_i\}$. If these $n$ exact ideal polygons does not cover all $\Omega$, we can add a few similar exact polygons where $\beta$ curves are outside of $\Omega$ to cover the missing parts in the middle of $\Omega$. Here, the main difference with the covering given in the lemma is that in this special covering all $\beta$ curves of $\D_i$ are outside of $\Omega$. 
\end{rmk}

\begin{rmk} [Exceptional Curves] \label{exception} Throughout the paper, we will not consider the infinite curves which are neither fat nor skinny. These are the curves which contains an ideal polygon $\Delta^\pm_i$ with $a(\Omega)=b(\Omega)$, but no skinny polygon at infinity. For example, Scherk curves described in section \ref{scherksec} are exceptional (See also Remark \ref{scherkrem}). After ignoring these exceptional curves, if an infinite curve is not fat, it must be skinny.
\end{rmk}

\section{Classification of Strongly Fillable Curves} \label{sfillsec}

In this section, we will prove our main result. In particular, we will study the role of $\Gamma^\pm$ on $\Gamma$ for strong fillability, and show that it completely determines strong fillability for a given infinite tall curve $\Gamma$. 

On the other hand, we will study the fillability question in the next section, and see that the fillability question and the strong fillability question are quite different for infinite curves. While $\Gamma^\pm$ completely determines strong fillability for a given infinite tall curve $\Gamma$, it is not very useful to detect fillability (See Section \ref{fillsec}).




Now, we will prove the main theorem. Recall that by curve, we mean a finite collection of disjoint Jordan curves. Let $\Gamma$ be an infinite curve in $\SI$. Let $\Gamma^\pm=\Gamma\cap (\overline{\BH^2}\times\{\pm\infty\})$ and $\wt{\Gamma}=\Gamma-(\Gamma^+\cup\Gamma^-)$ as before.

Note also our convention for {\em Exceptional Curves} (Remark \ref{exception}). In particular, in the following theorem, we will not consider the curves, which are neither fat nor skinny.

\begin{thm} \label{main} Let $\Gamma$ be a tame infinite curve in $\SI$. Then, $\Gamma$ is strongly fillable if and only if all of the following conditions satisfied:
	\begin{itemize}
		\item $int(\Gamma^\pm)$ is a collection of geodesics (possibly empty).
		\item $\Gamma$ is tall.	
		\item $\Gamma$ is nonoverlapping at the corner.
		\item $\Gamma$ is fat at infinity.
	\end{itemize}
\end{thm}

\begin{pf} We will divide the proof into two parts: In Step 1, we will deal with the "only if" part. In particular, we will show the nonexistence of area minimizing surfaces if $\Gamma$ does not satisfy one of the 4 conditions. In step 2, we will show the "if" part. In other words, we will show the existence of area minimizing surfaces if $\Gamma$ satisfy all the 4 conditions.

\vspace{.2cm}

\noindent {\bf Step 1:} If $\Gamma$ is strongly fillable, then $\Gamma$ must satisfy all of the 4 conditions.

\vspace{.2cm}

We naturally separate this step into 4 cases:

\vspace{.2cm}

\noindent {\bf Step 1a:} If $\Gamma$ is  strongly fillable,  then $\Gamma^\pm_g$ is a collection of geodesics.

\vspace{.2cm}

{\em Proof of Step 1a:} By Lemma \ref{geod}, for any fillable $\Gamma$, $\Gamma^\pm_g$ must be a collection of geodesics in $\caps$. \hfill $\Box$

\vspace{.2cm}

\noindent {\bf Step 1b:} If $\Gamma$ is strongly fillable, then $\Gamma$ must be tall.

\vspace{.2cm}

{\em Proof of Step 1b:} By Lemma \ref{finiteAPP}, for any strongly fillable $\Gamma$, $\Gamma$ must be tall, i.e. $S^1_\infty\times\BR-\wt{\Gamma}$ can be covered by tall rectangles. \hfill $\Box$

\vspace{.2cm}

\noindent {\bf Step 1c:} If $\Gamma$ is strongly fillable, then $\Gamma$ is nonoverlapping at the corner.

\vspace{.2cm}

{\em Proof of Step 1c:} Assume $\Sigma$ is an area minimizing surface in $\BHH$ with $\PI\Sigma=\Gamma$. As $\Gamma$ is a finite collection of Jordan curves in $\SI$, $\Gamma^\pm_c=\Gamma\cap (S^1_\infty \times \{\pm\infty\})$ is closed.

We claim $\Gamma^\pm_c$ does not contain any interval in $\CS$. Without loss of generality, assume $I=(p,q)\subset\Gamma^+_c$ where $I$ is an open interval in $S^1_\infty \times \{+\infty\}$. Let  $I'=[p',q']$ be a closed subinterval of $I$. As $\Gamma$ is a Jordan curve, for sufficiently large $c>0$, the rectangle $\widehat{\R}=I'\times[c,\infty)\subset \Si$ is disjoint from $\Gamma$.

For $t\in (\pi,\infty)$, let $\R_t= \partial (I'\times [c,c+t])$. In particular, if $h(t)=t$ is the height of the rectangle $\R_t$, then $h(t)\nearrow\infty$ as $t\nearrow \infty$. By Lemma \ref{infrec}, each $\R_t$ bounds a unique area minimizing surface $\T_t$ with $\PI \T_t=\R_t$ for $t\in(\pi,\infty)$. By \cite[Lemma 2.19]{Co}, as the asymptotic boundaries are disjoint ($\Gamma\cap \R_t=\emptyset$), then the area minimizing surfaces they bound are disjoint, i.e. $\Sigma\cap\T_t=\emptyset$.

Recall the infinite rectangles from Definition \ref{infrecdef}. As $h(t)\to \infty$ when $t\to \infty$, by \cite[Prop. 2.1]{ST}, $\T_t$ converges to an infinite rectangle $\widetilde{\T}$. In particular, here $t\to\infty$ corresponds to $d\searrow 1$ case in the proof of \cite[Prop. 2.1]{ST}. Note that $\PI \widetilde{\T}=\widetilde{\R}$ is a union of a pair of vertical line segments $\{p',q'\}\times[c,\infty]$, and a horizontal line segment $[p',q']\times\{c\}$ in $S^1_\infty\times\overline{\BR}$, and a geodesic segment $\gamma'$ in $\BH^2\times\{+\infty\}$ where $\PI\gamma'=\{p',q'\}$.

By assumption, $I'$ is in $\Gamma^+_c$. This implies $\PI\Sigma\supset I'$. However, by construction $\T_t\cap \Sigma=\emptyset$ for any $t\in (1,\infty)$. Again by construction $\widetilde{\T}\cap\Sigma\neq \emptyset$ as $\widetilde{\R}$ separates $\Gamma$ in $\SI$. Hence, as $\T_t\to \widetilde{\T}$ when $t\to \infty$, for sufficiently large $t_0$, $\T_{t_0}\cap\Sigma\neq \emptyset$. This is a contradiction which shows that $\Gamma^\pm_c$ does not contain any interval in $S^1_\infty \times \{\pm\infty\}$. Step 1c follows. \hfill $\Box$

\vspace{.2cm}

\noindent {\bf Step 1d:} If $\Gamma$ is strongly fillable, then $\Gamma$ is fat at infinity.

\vspace{.2cm}

{\em Proof of Step 1d:} Let $\Gamma$ be strongly fillable curve, i.e. $\PI \Sigma=\Gamma$ where $\Sigma$ is area minimizing surface in $\BHH$. Recall that as indicated in Remark \ref{exception}, we omit the infinite curves neither skinny nor fat. So, we will assume $\Gamma$ is skinny, and get a contradiction.

Assume that $\Gamma$ is skinny. Without loss of generality, assume $\Gamma^+=\gamma^1\cup..\cup\gamma^k$ induces an ideal $2k$-gon $\Delta$ where $a(\Delta)>b(\Delta)$. If $\Delta$ decomposes into more than one inscribed polygons, take $\Delta$ as the skinny polygon in the decomposition.

Let $\Sigma$ be an area minimizing surface in $\BHH$ with $\PI\Sigma=\Gamma$. Then, consider the sequence $S_n=\Sigma-n$ which is vertical translation down by $n$. By construction, the limit of the sequence $\{S_n\}$ is the collection of vertical geodesic planes $\Sigma^+=\Gamma^+_g\times\BR$ (Lemma \ref{geod}). Furthermore, as the limit of area minimizing surfaces is area minimizing, $\Sigma^+$ is also area minimizing.

Let $\partial \Delta=\gamma^1\cup\beta^1\cup ... \cup\gamma^k\cup \beta^k$. By assumption $a(\Delta)>b(\Delta)$. Here, $a(\Delta)$ corresponds to total "length" of $\Gamma^+=\gamma^1\cup..\cup\gamma^k$, and $b(\Delta)$ corresponds to the total "length" of the remaining geodesics in $\partial \Delta$, i.e. $\partial \Delta^+-\Gamma^+=\beta^1\cup...\beta^k$ (See Section \ref{scherksec}).

Recall that $B_m$ be the disk of radius $m$ and center $O$ in $\BH^2$, and $\B_m=B_m\times[-m,m]$ is the solid cylinder in $\BHH$. Consider $\Sigma^+_m=\B_m\cap \Sigma^+$. We claim that for sufficiently large $m$, $\Sigma^+_m$ is not an area minimizing surface. Let $\eta_m=\partial \Sigma^+_m$ be the collection of disjoint $k$ Jordan curves in $\partial \B_m$.

For $1\leq i\leq k$, let $\gamma^i_m=\gamma^i\cap B_m$. In other words, $\gamma^i_m$ is a finite arc segment in the infinite geodesic $\gamma^i$. Then, $\Sigma^+_m=\bigcup \gamma^i_m\times[-m,m]$ is a collection of $k$ vertical geodesic surfaces in $\B_m$ by construction. Let $\partial\gamma^i_m=\{p^{2i-1}_m,p^{2i}_m\}$ be the endpoints of $\gamma^i_m$. Then, we have $2k$ points $\V_m=\{p^1_m,p^2_m,...,p^{2k}_m\}$ in $\partial B_m$. Let $\beta^i_m$ be the geodesic in $B_m$ connecting $p^{2i}_m$ and $p^{2i+1}_m$. Hence, $\gamma^i_m$ and $\beta^i_m$ curves defines a polygon $\Delta_m$ in $B_m$, i.e. $\partial\Delta_m=\gamma^1_m\cup\beta^1_m\cup...\gamma^k_m\cup\beta^k_m$.

Now, let $\Pi^i_m=\beta^i_m\times [-m,m]$ be a vertical geodesic surface in $\B_m$. Let $\Delta^+_m=\Delta_m\times\{m\}$ and $\Delta^-_m=\Delta_m\times\{-m\}$, i.e. $\Delta^\pm_m\subset \partial\B_m$. Then, define a surface $S_m=\bigcup_{i=1}^k \Pi^i_m\cup\Delta^+_m\cup\Delta^-_m$ in $\B_m$. $S_m$ is topologically a sphere with $k$ holes. Furthermore, $\partial S_m=\partial \Sigma_m=\eta_m$.

We claim that the area of $S_m$ is less than the area of $\Sigma_m$ for sufficiently large $m$. Let $\|.\|$ and $|.|$ represent the area and the length respectively. Since $\Delta$ is an ideal $2k$-gon, then $\|\Delta\|=2(k-1)\pi$. By construction, $\Delta_m^\pm\subset \Delta$ for any $m$. Hence, $\|\Delta^\pm_m\|<2(k-1)\pi$.

Now, by assumption $a(\Delta)>b(\Delta)$. By the definition of $a(.)$ and $b(.)$, this implies that $a_m=\Sigma_{i=1}^k|\gamma^i_m|>\Sigma_{i=1}^k|\beta^i_m|=b_m$ for sufficiently large $m$. Let $c_m=a_m-b_m$ for $m$ large, and let $c=a(\Delta)-b(\Delta)>0$. Then, $c_m\nearrow c$ as $m\to\infty$, and hence $c_m>0$ for sufficiently large $m$. Now, $\|\Sigma^+_m\|=2m.a_m$ and $\|S_m\|<2m.b_m+4(k-1)\pi$. Since $c_m\nearrow c$, for sufficiently large $m$, $\|\Sigma^+_m\|>\|S_m\|$. This proves that $\Sigma^+_m$ is not an area minimizing surface, and give a contradiction. Step 1d follows. \hfill $\Box$

This finishes the proof of Step 1. \hfill $\Box$

\

\noindent {\bf Step 2:} If $\Gamma$ satisfy all of the 4 conditions, then $\Gamma$ is strongly fillable.

\vspace{.2cm}

{\em Proof of Step 2:} We will show that there exists an area minimizing surface $\Sigma$ in $\BHH$ with $\PI\Sigma=\Gamma$. Naturally, we construct a sequence of compact area minimizing surfaces $\Sigma_n$ with $\partial \Sigma_n\to \Gamma$. Our aim is to take the limit of $\{\Sigma_n\}$, and to show that the limit area minimizing surface $\Sigma$ in $\BHH$ has the asymptotic boundary $\Gamma$. However, as indicated in \cite{Co}, the sequence might escape to infinity. Then, we might end up with an empty limit, or a nonempty limit $\Sigma$ with $\PI \Sigma\subset \Gamma$ but $\PI \Sigma\neq \Gamma$. Hence, to prevent the sequence escape to infinity, we first construct a barrier $\N$ near $\SI-\Gamma$ (See also \cite[Proposition 4.1]{KM}).


\vspace{.2cm}

\noindent {\bf Step 2a:} The construction of the barrier $\N$ near infinity.

\vspace{.2cm}

\noindent {\em Outline:} The barrier $\N$ can be considered as a neighborhood of $\Gamma^c=\SI-\Gamma$ in $\overline{\BHH}$. By using this barrier $\N$ near the asymptotic boundary, we will define a mean convex domain $\Omega=\N^c$ such that $\PI \Omega=\Gamma$ and $\Sigma_n\subset \Omega$. In particular, we want to keep the sequence $\{\Sigma_n\}$ away from $\SI-\Gamma$ in order to prevent $\{\Sigma_n\}$ {\em escape to infinity}, i.e. $\lim\Sigma_n=\Sigma\subset \Omega$ and $\PI\Sigma=\PI \Omega=\Gamma$. Hence, the condition $\N\cap\Sigma_n=\emptyset$ makes sure this, and $\N$ would act as a barrier between $\{\Sigma_n\}$ and $\SI-\Gamma$. In \cite{Co}, we constructed such a barrier at infinity for finite curves. Now, we construct a similar barrier for infinite curves. See also the "summary of the barrier construction" at the end of the proof.

We will use the notation above. To construct the barrier near the cylinder $\Si$, we can use the tall rectangles again as in \cite{Co}. Hence, the main problem is the constructing a barrier near the caps at infinity $\caps$ to prevent the sequence from escaping to infinity.

$\Gamma$ is a finite collection of disjoint Jordan curves in $\SI$ with $\Gamma=\Gamma^+\cup\Gamma^-\cup\wt{\Gamma}$. Without loss of generality, we assume either $\Gamma_g^+$ or $\Gamma_g^-$ is nonempty, and $\Gamma_g^+=\gamma_1^+\cup ...\cup\gamma^+_n$, and $\Gamma_g^-=\gamma_1^-\cup ...\cup\gamma^-_m$. We will deal with the trivial cases $m=0,1$ or $n=0,1$ at the end. So, we will assume $m,n>1$. 

Define $\Delta^\pm$, $\tau_i^+$ and $\tau_j^-$ as in Section \ref{fatsec} where $1\leq i\leq 2n$ and $1\leq j\leq 2m$. In particular, $\Delta^+$ is an ideal $2n$-gon and $\Delta^-$ is an ideal $2m$-gon in $\caps$ defined by $\{\gamma_i^\pm\}$, and $\tau_i^\pm$ represent all infinite geodesics of $\partial \Delta^\pm$. 

Furthermore,  let $\wh{\Delta}^\pm$ be the convex hull of $\Gamma^\pm=\Gamma^\pm_g\cup\Gamma^\pm_c$ in $\overline{\BH^2}\times\{\pm\infty\}$. In particular, $\Delta^\pm\subseteq\wh{\Delta}^\pm$, and if $\Gamma^\pm_c=\emptyset$, then $\wh{\Delta}^\pm=\Delta^\pm$. 

The barrier $\N$ consists of three major blocks: 
\begin{itemize}
	\item Infinite side barriers: Covering $\caps-\wh{\Delta}^\pm$
	\item Scherk barriers at infinity: Covering $\wh{\Delta}^\pm-\Gamma^\pm_g$
	\item Tall rectangles: Covering $\Si-\wt{\Gamma}$
	
\end{itemize}



\noindent {\em Infinite side barriers:} In this part, we want to cover the complement of $\wh{\Delta}^\pm$ in $\overline{\BH^2}\times\{\pm\infty\}$. Let $\partial \wh{\Delta}^+=\bigcup_{i=1}^{n'}\mu_i^+$ where $n'=2n+ \#(\Gamma^+_c)$. For each $\mu^+_i$, we define an infinite side barrier $\wh{\T}_i^+$ as follows. Let $\PI\mu_i^+=\{p_i^+,p_{i+1}^+\}$ where $\{p_i^+\}$ are circularly ordered as before. Each $\mu^+_i$ in $\partial \wh{\Delta}^+$, defines a lens shaped region $\U^+_i$ in $\BH^2\times\{+\infty\}$ where $\wh{\Delta}^+\cap\overline{\U^+_i}=\mu_i^+$ and $\BH^2\times\{+\infty\}-\wh{\Delta}^+=\bigcup_1^{n'} \U_i^+$, i.e. $\partial \overline{\U^+_i}=\mu_i^+$ and $\PI\overline{\U^+_i}=[p_i^+,p_{i+1}^+]$ (See Figure \ref{scherkbarrier}-left). As $\U_i^+ \subset \SI-\Gamma$, we would like to cover $\U_i^+$ with infinite rectangles.

Let $p_i^+=(\theta_i,+\infty)$ for $1\leq i\leq n'$ where $\theta_i\in [0,2\pi)$. If there exists an infinite rectangle $\R_i^+$ in $\SI$ (Definition \ref{infrecdef}) such that $\R_i^+\cap \Gamma = \mu_i^+$, then let $\T_i^+$ be the unique area minimizing surface $\R_i^+$ bounds in $\BHH$ by Lemma \ref{infrec}. Let $\wh{\U}_i^+$ be the open domain separated by $\T^+_i$ from $\BHH$ with $\PI \wh{\U}^+_i\supset \U_i^+$ (See Figure \ref{sidecover}-left).

If there is no such $\R_i^+$, we will cover $\U_i^+$ as follows. Recall that $\PI \overline{\U_i^+}=[p_i^+,p_{i+1}^+]=[\theta_i,\theta_{i+1}]\times\{+\infty\}$ in $S^1_\infty\times\{+\infty\}$. Recall that by construction,  $((\theta_i,\theta_{i+1})\times\{+\infty\})\cap\Gamma=\emptyset$. Fix $\e_o>0$ sufficiently small. For any $t\in(0,\e_0)$, let $N_t=\sup\{N\mid \Gamma\cap [\theta_i+t,\theta_{i+1}-t]\times[N,\infty)\neq\emptyset\}$. As $((\theta_i,\theta_{i+1})\times\{+\infty\})\cap\Gamma=\emptyset$, for any $t\in(0,\e_0)$, we have $N_t<\infty$. Define $\R^+_{it}$ to be the infinite rectangle containing $[\theta_i+t,\theta_{i+1}-t]\times[N_t+1,\infty)$. By construction, $\R^+_{it}\cap\Gamma=\emptyset$ (See Figure \ref{sidecover}-right).

\begin{figure}[t]
\begin{center}
$\begin{array}{c@{\hspace{.4in}}c}

\relabelbox  {\epsfysize=2in \epsfbox{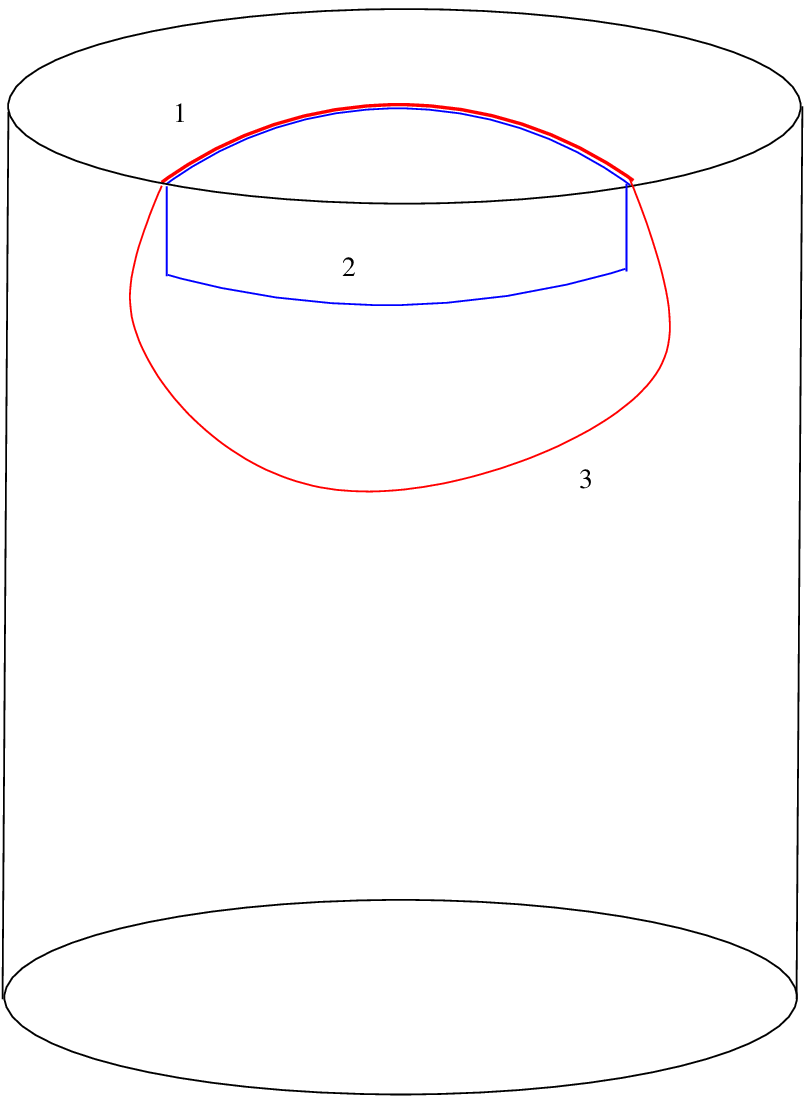}}
\relabel{1}{\scriptsize $\mu^+_i$} \relabel{2}{\scriptsize $\R_i^+$} \relabel{3}{\scriptsize $\Gamma$}
 \endrelabelbox &

\relabelbox  {\epsfysize=2in \epsfbox{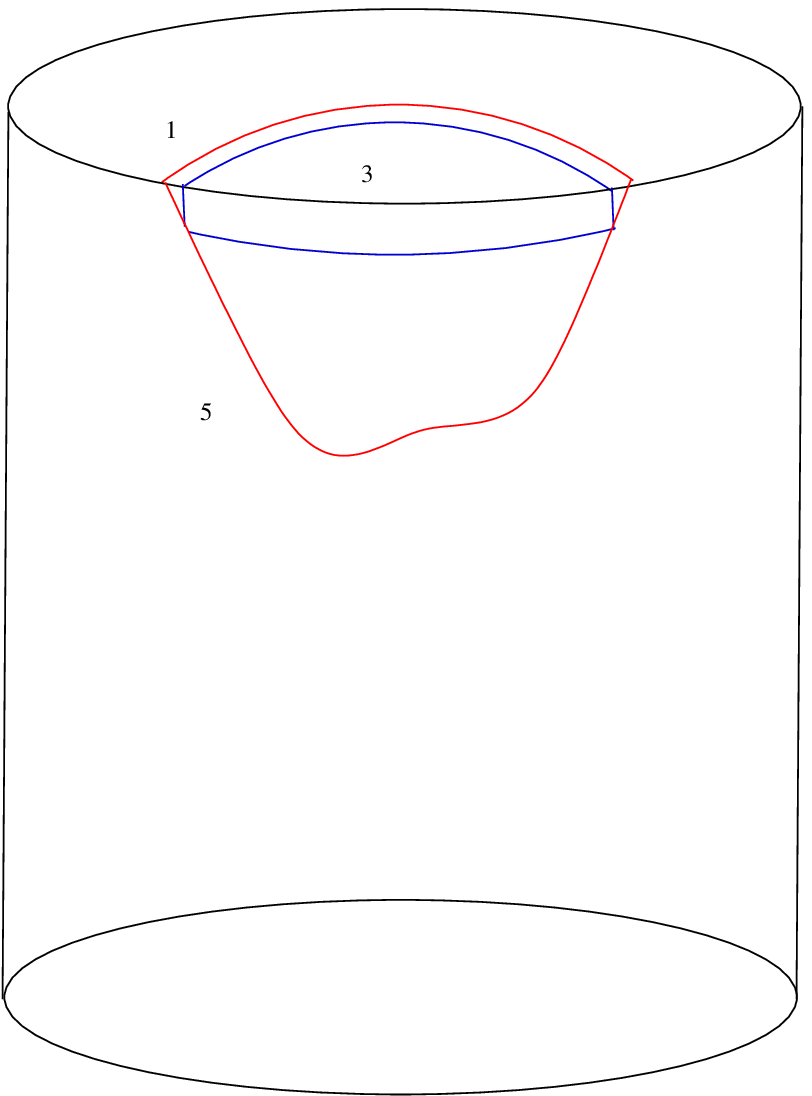}}
\relabel{1}{\tiny $\mu^+_i$}   \relabel{3}{\tiny $\R^+_{it}$}
 \relabel{5}{\scriptsize $\Gamma$}
\endrelabelbox \\
\end{array}$

\end{center}

\caption{ \label{sidecover} \footnotesize In the left, we have the trivial case, where we can cover outside of $\mu^+_i$ just by one infinite rectangle $\R_i^+$ as $\wt{\Gamma}$ curves towards "outside". In the right, we need to use a family of infinite rectangles $\{\R_{it}^+\}$ for covering as $\wt{\Gamma}$ curves towards "inside".}
\end{figure}

Now, let $\T^+_{it}$ be the unique area minimizing surface bounding $\R^+_{it}$. Let $\wh{\U}^+_{it}$ be the open component of $\BHH-\T^+_{it}$ where $\PI \wh{\U}^+_{it}$ contains the arc $(\theta_i+t, \theta_{i+1}-t)\times\{+\infty\}$ in the upper corner circle of $\SI$. Then define $\wh{\U}^+_i=\bigcup_t \wh{\U}^+_{it}$. Notice that $\PI \wh{\U}^+_i\cap \BH^2\times\{+\infty\}=\U^+_i$. Furthermore, $\PI \wh{\U}^+_i\supset ((\theta_i,\theta_{i+1})\times\{+\infty\})$. Similarly define $\wh{\U}^-_j$ for $1\leq j\leq m'$. We will call $\wh{\U}^\pm_i$ as infinite side barrier. Notice that for any $\mu^\pm_i$, we have an infinite side barrier $\wh{\U}^\pm_i$ such that $\Gamma\cap \PI \wh{\U}^\pm_i=\emptyset$. Hence, infinite side barriers cover outside of $\wh{\Delta}^\pm$ in $\caps$.

\vspace{.2cm}

\noindent {\em Scherk barriers at infinity:} Now, we are going to construct the second major block of our barrier $\N$. After covering outside of $\wh{\Delta}^\pm$ by infinite side barriers, we want to cover the inside of $\wh{\Delta}^\pm- \Gamma^\pm_g$ to construct $\N$ so that $\PI\N^c=\Gamma$. We will construct this in two steps. First, we will cover inside of $\Delta^\pm$, i.e. convex hull of $\Gamma^\pm_g$. Then, we will cover the remaining part $\wh{\Delta}^\pm-\Delta^\pm$.

First, consider $\Delta^\pm$. By Section \ref{fatsec}, $\Delta^\pm$ decomposes into inscribed polygons by $\Gamma^\pm$, i.e. $\Delta^+=\Delta^+_1\cup...\Delta^+_{c^+}$ and $\Delta^-=\Delta^-_1\cup...\Delta^-_{c^-}$ where $c^\pm\geq 1$. Since $\Gamma$ is fat at infinity,  $\Delta^\pm_i$ is a fat polygon for any $i$, and it has an exact covering by Lemma \ref{fatlem} and Remark \ref{specialcover}. Let $\{\D^+_{ik}\}$ be the special exact covering described in Remark \ref{specialcover} such that for any $1\leq i\leq c^+$, $\Delta^+_i \subset\bigcup_{k=1}^{n_i} \D^+_{ik}$. Furthermore, all $\beta$ curves of $\D^+_{ik}$ is outside of $\Delta_i^+$ (See Figure \ref{fatcovering}-right). 

Fix $\Delta^+_{i_o}$. Notice that $\Delta^+_{i_o}$ is an ideal $2n_{i_o}$-gon for some $n_{i_o}\leq n$, and $\Delta^+_{i_o} \subset\bigcup_{k=1}^{n_{i_o}} \D^+_{i_o k}$. By construction, $\partial \Delta^+_{i_o}=\gamma_{j_1}\cup\tau_{j_1}\cup...\cup\gamma_{j_{ni_o}}\cup\tau_{j_{ni_o}}$. Recall that $\alpha$ curves of $\Delta^+_{i_o}$ are the $n_{i_o}$ geodesics $\gamma_{j_1},..,\gamma_{j_{ni_o}}$ in $\Gamma^+ \cap \partial \Delta^+_{i_o}$, while $\beta$ curves of $\Delta^+_{i_o}$ are the remaining $n_{i_o}$ geodesics $\tau_{j_1},..,\tau_{j_{ni_o}}$. Recall that in our special covering only $\alpha$ curves of $\D^+_{i_o k}$ intersects $\Delta^+_{i_o}$, while $\beta$ curves of $\D^+_{i_o k}$ stays completely outside of $\Delta^+_{i_o}$ (See Figure \ref{fatcovering}-right). In particular, $\bigcup_{i,k} int(\D_{ik}^+)\supset \Delta^+-\Gamma^+$ and $\bigcup_{i,k} int(\D_{ik}^+)\cap\Gamma^+=\emptyset$.

\begin{figure}[h]
\begin{center}
$\begin{array}{c@{\hspace{.4in}}c}

\relabelbox  {\epsfysize=2in \epsfbox{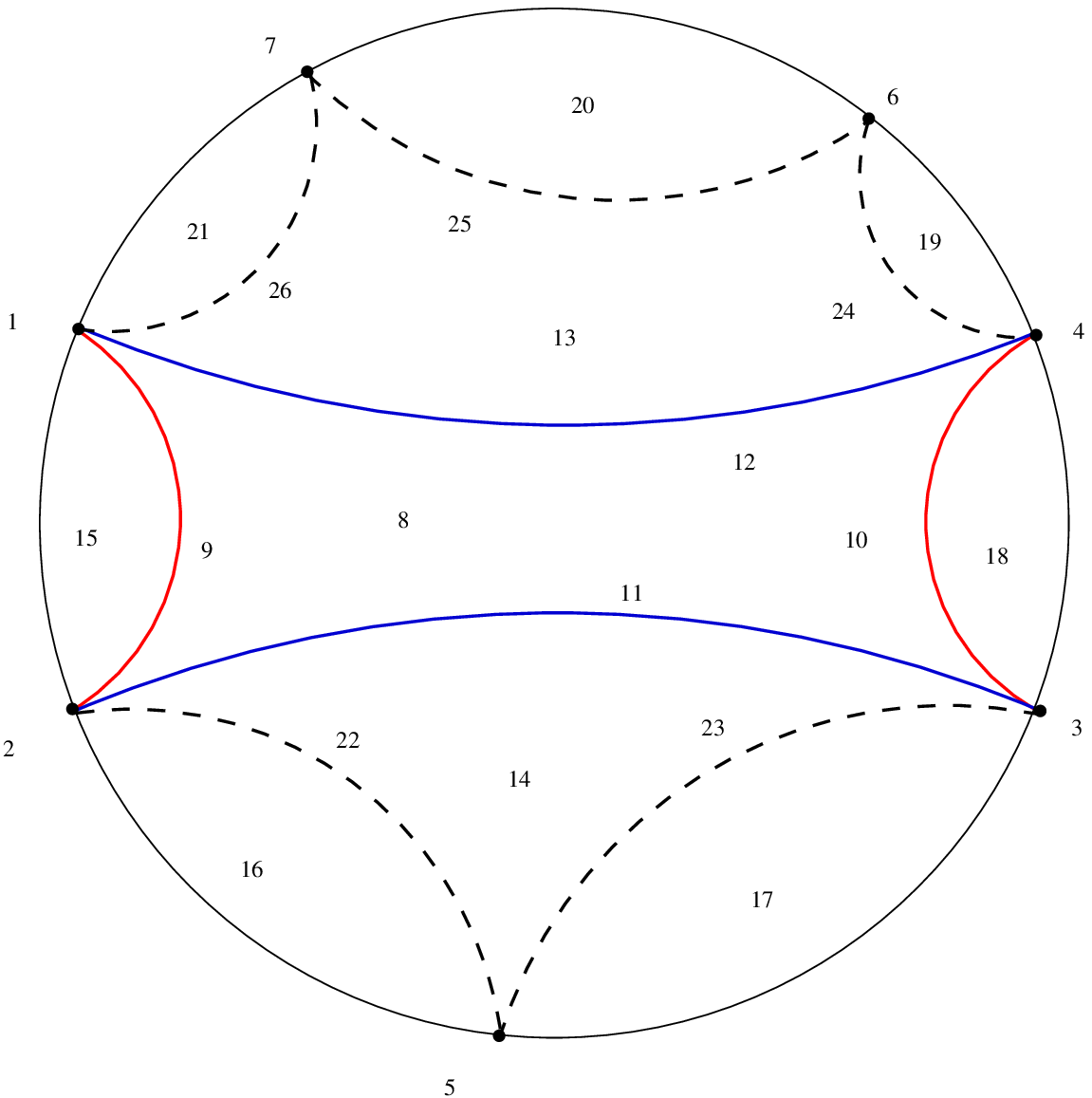}}
\relabel{1}{\footnotesize $p_1$} \relabel{2}{\footnotesize $p_2$ }  \relabel{3}{\footnotesize $p_3$} \relabel{4}{\footnotesize $p_4$} \relabel{5}{\footnotesize $p_5$} \relabel{6}{\footnotesize $p_6$ }  \relabel{7}{\footnotesize $p_7$} \relabel{8}{\footnotesize $\Delta^+$}

\relabel{9}{\footnotesize $\alpha_1$} \relabel{10}{\footnotesize $\alpha_2$ }  \relabel{11}{\footnotesize $\beta_1$} \relabel{12}{\footnotesize $\beta_2$}

\relabel{13}{\footnotesize $\V_4$} \relabel{14}{\footnotesize $\V_2$} 

\relabel{15}{\scriptsize $\U_1$} \relabel{16}{\scriptsize $\U_2$} \relabel{17}{\scriptsize $\U_3$}
\relabel{18}{\scriptsize $\U_4$} \relabel{19}{\scriptsize $\U_5$} \relabel{20}{\scriptsize $\U_6$} \relabel{21}{\scriptsize $\U_7$}

\relabel{22}{\footnotesize $\mu_2$} \relabel{23}{\footnotesize $\mu_3$} \relabel{24}{\footnotesize $\mu_5$}
\relabel{25}{\footnotesize $\mu_6$} \relabel{26}{\footnotesize $\mu_7$} 

\endrelabelbox &

\relabelbox  {\epsfysize=2in \epsfbox{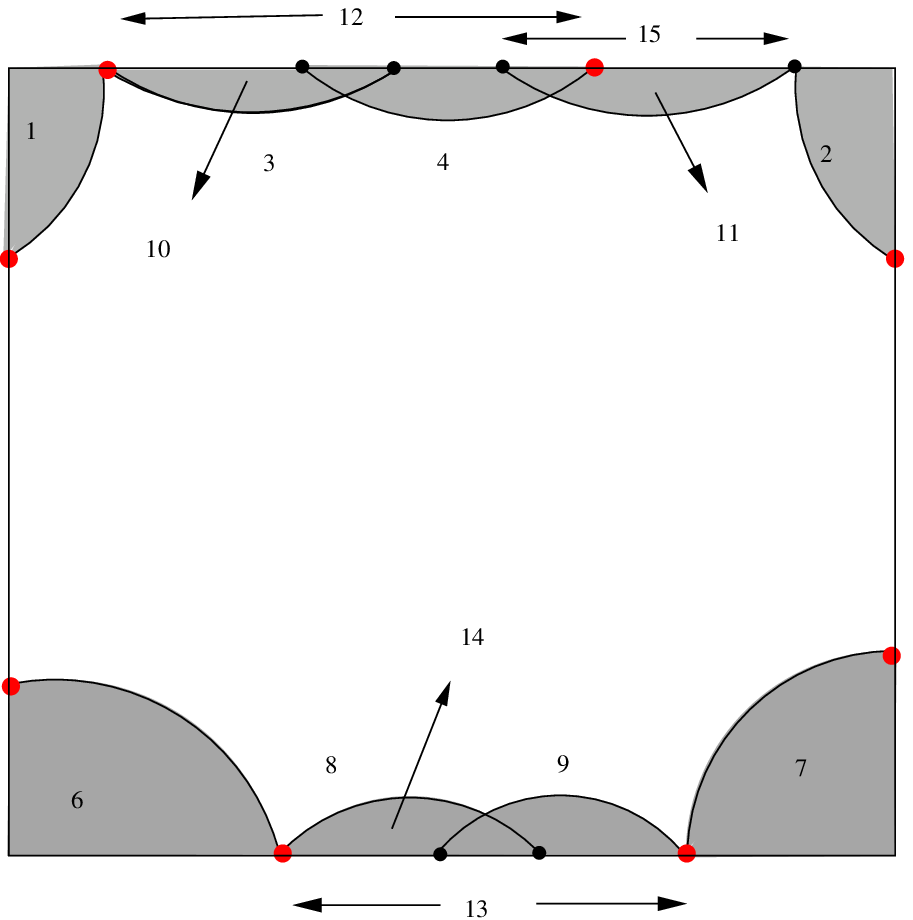}}
\relabel{1}{\tiny $\wh{\U}^+_1$} \relabel{2}{\tiny $\wh{\U}^+_3$}
\relabel{3}{\tiny $S^t_{i 1}$}
\relabel{4}{\tiny $S^t_{i 2}$}
\relabel{6}{\tiny $\wh{\U}^-_1$} \relabel{7}{\tiny $\wh{\U}^-_3$}
\relabel{8}{\tiny $S^{-t}_{j 1}$} \relabel{9}{\tiny $S^{-t}_{j 2}$}
\relabel{10}{\tiny $\Y^+_{i 1}$} \relabel{11}{\tiny $\W_{j1}^+$}
\relabel{12}{\tiny $\Delta_i^+$} \relabel{13}{\tiny $\Delta_k^-$} 

\relabel{14}{\tiny $\Y_{k1}^+$} \relabel{15}{\tiny $\V_j^+$} 

\endrelabelbox \\
\end{array}$

\end{center}

\caption{ \label{scherkbarrier} \footnotesize In the left, we have $\Gamma^+_g=\alpha_1\cup\alpha_2$ and $\Gamma^+_c=\{p_5,p_6,p_7\}$. This gives a decomposition as $\BH^2\times\{+\infty\}=\wh{\Delta}^+\bigcup\U_i$ and $\wh{\Delta}^+=\Delta^+\bigcup\W_j$. In the right, red points represent $\Gamma$, and a cross-section of $\N^+$ and $\N^-$ are pictured from the side.}
\end{figure}

Now, we will cover $\wh{\Delta}^\pm- \Delta^\pm$ by exact polygons. Recall that $\Delta^\pm$ is the convex hull of $\Gamma^\pm_g$ in $\caps$. We abuse the notation by using a different indexing for vertices $\{p_i^+\}$ of $\Delta^+$ and vertices $\{p_i^+\}$ of $\wh{\Delta}^+$ in previous section when defining $\{\mu_i^+\}$. Let $\PI\Gamma^+_g=\{p_1^+,...,p_{2n}^+\}$ be circularly ordered. Then, $\partial \Delta^+= \tau_1^+\cup\tau_2^+\cup...\cup\tau_{2n}^+$ where $\tau_i^+=\overline{p_i p}_{i+1}$ (In Figure \ref{scherkbarrier}-left, $\tau_1=\alpha_1$, $\tau_2=\beta_1$, $\tau_3=\alpha_2$, and $\tau_4=\beta_2$). For each $\tau_i^+$, consider the component $\V_i^+$ in $\wh{\Delta}^+- \Delta^+$ such that $\partial \overline{\V_i^+}\supset \tau_i^+$. Notice that if $\Gamma^+_c\cap (p_i,p_{i+1})=\emptyset$, then $\V_i^+=\emptyset$. If $\Gamma^+_c\cap (p_{i_o},p_{i_o+1})=\{q_1,q_2,...q_{e_o}\}$, then $\V_{i_o}^+$ is the convex hull of $\{p_{i_o},p_{i_o+1},q_1,...,q_{e_o}\}$ in $\BH^2\times\{+\infty\}$.

Now, fix $1\leq i_o\leq 2n$, and assume $\Gamma^+_c\cap (p_{i_o},p_{i_o+1})=\{q_1,q_2,...q_{e_o}\}$. Define a covering of $\V_{i_o}$  by exact $(2{e_o}+2)$-gons $\{\D_{i_o k}^+\}$ such that for any $k$, $\tau_{i_o}^+$ is an $\alpha$-curve of $\D_{i_o k}^+$, and all $\beta$-curves of $\D_{i_o k}^+$ are outside of $\wh{\Delta}^+$. We can obtain such a covering as follows: Let $(p_{i_o},p_{i_o+1})-\{q_1,...q_{e_o}\}=I_0\cup I_1\cup...\cup I_{e_o}$ where $I_j$ is an open interval in $S^1_\infty$. Let $\{p_{i_o},p_{i_o+1}\}$ be the first two vertices of the exact $(2{e_o}+2)$-gon $\D_{i_o k}^+$. Then, place remaining $2{e_o}$ vertices of $\D_{i_o k}^+$ such that let one vertex $x_1$ be in $I_0$ and let one vertex $x_{2{e_o}}$ be in $I_{e_o}$. Then place two vertices $\{x_{2j},x_{2j+1}\}$  in $I_j$ for each $1\leq j \leq {e_o}-1$. Notice that $q_j\in (x_{2j-1},x_{2j})$ for any $1\leq j \leq {e_o}$. Then, the $\alpha$ curves of $\D_{i_o k}^+$ would be $\tau_{i_o}^+$ and the geodesics $\overline{x_{2j-1} x}_{2j}$ for $1\leq j \leq {e_o}$. Similarly, the $\beta$ curves of $\D_{i_o k}^+$ would be the geodesics $\overline{x_{2j} x}_{2j+1}$ for $1\leq j \leq {e_o-1}$,  $\overline{p_i x}_1$, and $\overline{x_{2e_o} p}_{i+1}$. Since we can choose the vertices $x_{2j-1}$ and $x_{2j}$ as far or close as we want to the vertex $q_j$, we can make  $(2{e_o}+2)$-gons $\{\D_{i_o k}^+\}$ both exact and covering $\V_{i_o}^+$. Hence, this gives a covering of $\wh{\Delta}^\pm- \Delta^\pm$ by exact polygons. In particular, $\bigcup_{i,k} int(\D_{ik}^+)\supset \wh{\Delta}^+-\Delta^+$ and $\bigcup_{i,k} int(\D_{ik}^+)\cap\Gamma^+=\emptyset$.

So far, we constructed a family of exact polygons covering $\wh{\Delta}^\pm -\Gamma^\pm$. We will use these exact polygons to build Scherk barriers near infinity. We will start with the exact polygons covering $\wh{\Delta}^\pm- \Delta^\pm$. Let $\D^+_{i_o k}$ be such an exact polygon used in the covering $\V_{i_o}^+$ in $\wh{\Delta}^+ - \Delta^+$. Consider the Scherk graph $S_{i_o k}$ (Lemma \ref{scherk}) over $\D^+_{i_o k}$ where it takes $+\infty$ value on the $\alpha$-curves of $\D^+_{i_o k}$, and $-\infty$ value on the $\beta$-curves. Recall that we have the infinite side barriers $\wh{\U}^+_{i}$ for each $\mu^+_i\subset \partial \wh{\Delta}^+$. Let $S^t_{i_o k}=S_{i_o k}+t$ be the vertical translation of $S_{i_o k}$ by $t$ units up. 

Notice that as $t\to +\infty$, $S^t_{i_o k}\to \wh{\beta}_{i_o k}\times \BR$ where $\wh{\beta}_{i_0 k}$ is the collection of $\beta$-curves of $\D^+_{i_o k}$. Since all $\beta$ curves of $\D^+_{i_o k}$ are outside of $\wh{\Delta}^+$, for sufficiently large $t_o$, the Scherk graph $S^{t_o}_{i_o k}$ intersects $\partial \mathbf{cl}(\wh{\U}^+_{j_i})$ in an infinite arc $\sigma_{i_o k}^{j_i}$ with $\PI\sigma_{i_o k}^{j_i}\subset \mu_{j_i}$ for some $j_i$ with $\D^+_{i_o k}\cap\mu^+_{j_i}\neq \emptyset$. Here $\mathbf{cl}(.)$ means closure.

Let $\wh{S}_{i_o k}=S^{t_o}_{i_o k}- \bigcup_j \wh{\U}_{j}^+$. Intuitively, $\wh{S}_{i_o k}$ is a subsurface of the Scherk surface $S^{t_o}_{i_o k}$ where the parts going to $-\infty$ are cut out by the infinite curves $\sigma_{i_o k}^{j_i}$. Hence, $\PI \wh{S}_{i_o k} \subset \BH^2\times\{+\infty\}$, and $\partial \wh{S}_{i_o k}= \bigcup_j\sigma_{i_o k}^{j}$. Let $\W_{i_o k}$ be the open connected region over $\wh{S}_{i_o k}$ where $\PI \W_{i_o k}\subset \BH^2\times\{+\infty\}$ and $\partial \overline{\W}_{i_o k} \subset \wh{S}_{i_o k}\bigcup_j \partial \mathbf{cl}(\wh{\U}^+_{j})$. Hence, we define a barrier $\W^+_{i_o k}$ near the upper cap for each exact polygon $\D^+_{i_o k}$ in the covering of $\V_{i_o}$ which is a component of $\wh{\Delta}^\pm -\Delta^\pm$. 
By the construction and its shape, we call $\W^+_{i_o k}$  Scherk barrier at infinity. Let $\W_{i_0}^+= \bigcup_k \W^+_{i_o k}$ be the Scherk barrier corresponding to $\V_{i_o}^+$. In particular, $\PI \W_{i_0}^+ \subset \BH^2\times \{+\infty\}$ and $\V_{i_o}^+ \subset \PI \W_{i_0}^+$.

Now, we will finish the construction with a similar process for the exact polygons $\Delta^\pm-\Gamma^\pm$. Again, let $\D^+_{i_1 k}$ be such an exact polygon used in the covering of the fat polygon $\Delta_{i_1}^+$ in $\Delta^+$. Consider the Scherk graph $S_{i_1 k}$ over $\D^+_{i_1 k}$ where it takes $+\infty$ value on the $\alpha$-curves of $\D^+_{i_1 k}$, and $-\infty$ value on the $\beta$-curves. Recall that the $\alpha$-curves of $\D^+_{i_1 k}$ are $\partial\D^+_{i_1 k}\cap \Gamma^+$, and $\beta$-curves are the remaining curves in $\partial\D^+_{i_1 k}$, which are outside of $\Delta_{i_1}^+$ by construction.

Now, similar to above construction, we translate $S_{i_1 k}$ sufficiently up so that $S^{t_1}_{i_1 k}$ intersects $\bigcup_j\wh{\U}_j^+\bigcup_i\W_i^+$ in an infinite arc $\sigma_{i_1 k}^{j_i}$ with $\PI\sigma_{i_1 k}^{j_i}\subset \tau_{j_i}\subset \partial \Delta^+_{i_1}$. Notice that if $\V_{j_i}$ is nontrivial, $S^{t_1}_{i_1 k}$ intersects the corresponding Scherk barrier $\W_{j_i}$. Similarly, if $\V_{j_i}=\emptyset$, $S^{t_1}_{i_1 k}$ intersects the corresponding infinite side barrier $\wh{U}_{j_i}^+$. Again, we get a subsurface $\wh{S}_{i_1 k}$ of the Scherk surface $S^{t_1}_{i_1 k}$ such that $\PI \wh{S}_{i_1 k} \subset \BH^2\times\{+\infty\}$, and $\partial \wh{S}_{i_1 k}= \bigcup_j\sigma_{i_1 k}^{j}$. Similarly, for any exact polygon $\D^+_{i_1 k}$, we define an open connected region above $\wh{S}_{i_1 k}$, say $\Y^+_{i_1 k}$. Then, we have $\Y^+_{i_1}=\bigcup \Y^+_{i_1 k}$ is the Scherk barrier corresponding to the fat polygon $\Delta_{i_1}^+$. Notice that $\PI \Y^+_{i_1}\subset \BH^2\times \{+\infty\}$ and $\Delta_{i_1}^+-\Gamma^+ \subset \PI \Y^+_{i_1}$.

Now, we complete the barrier construction near the caps at infinity: Define $\N^+=\bigcup_i \W_i^+ \bigcup_j \Y_j^+ \bigcup_k \wh{\U}_k^+$. In particular, we see that $\PI\N^+\cap \Gamma= \emptyset$, and furthermore, $\PI\N^+ \cap \overline{\BH^2}\times\{+\infty\}=\overline{\BH^2}\times\{+\infty\}-\Gamma^+$. This means $\N^+$ can be considered as a barrier for $\BH^2\times\{+\infty\}-\Gamma^+$ as well as its neighborhood in $\overline{\BHH}$. Define $\N^-$ similarly. See Figure \ref{scherkbarrier}-right.

\vspace{.2cm}

\noindent {\em Summary of the Barrier Construction:} Given $\Gamma^+=\Gamma^+_g\cup \Gamma^+_c$, we naturally define following regions: $\wh{\Delta}^+=CH(\Gamma^+)$ and $\Delta^+=CH(\Gamma^+_g)$ where $CH(.)$ represents the convex hull. This gives us the decomposition of $\BH^2\times\{+\infty\}$ into following domains (See Figure \ref{scherkbarrier}-left):

\begin{itemize}	
	\item $\U_k^+$: Lens shaped components of $\BH^2\times\{+\infty\}-\wh{\Delta}^+$.
	\item $\V_j^+$: Ideal polygons of $\wh{\Delta}^+-\Delta^+$.
	\item $\Delta_i^+$: Fat polygons in decomposition of $\Delta^+$ by $\Gamma^+_g$.	
\end{itemize}

For each of these domains in $\BH^2\times\{+\infty\}$, we defined the following corresponding neighborhoods in $\BHH$ (See Figure \ref{scherkbarrier}-right).

\begin{itemize}
	\item $\wh{\U}_k^+$: Infinite side barrier $\wh{\U}_k$ such that $\PI \wh{\U}_k^+ \supset \U_k^+$.
	\item $\W_j^+$: Scherk Barrier $\W_j^+$ such that $\PI \W_j^+\supset \V_j^+$.
	\item $\Y_i^+$: Scherk Barrier $\Y_i^+$ such that $\PI \Y_i^+\supset \Delta_i^+$.
\end{itemize}

Final piece of the barrier $\N$ is the one covering the finite part $\wt{\Gamma}$. We dealt with this case in \cite[Theorem 2.13]{Co} where we solve the asymptotic Plateau problem for finite curves. As $\Gamma$ is tall, we can cover $\Si - \wt{\Gamma}$ with tall rectangles $\R_k$ in $\Si$, i.e. $\Si - \wt{\Gamma}=\bigcup_k\R_k$. Let $\T_k$ be the unique area minimizing surface in $\BHH$ with $\PI \T_k=\partial\R_k$, and $\wh{\U}_k$ be the component of $\BHH-\T_k$ with $\PI\wh{\U}_k=int(\R_k)$. Then define $\wt{\N}=\bigcup_k\wh{\U}_k$, the third major block of $\N$. Notice that $\PI \wt{\N}=\Si - \wt{\Gamma}$ by construction.

Now, define $\N=\N^+\cup\N^-\cup\wt{\N}$. Then, $\N$ is an open region in $\BHH$ with $\PI\N=\SI-\Gamma$.

\vspace{.2cm}

\noindent {\em Mean Convexity:} Notice that $\Omega=\BHH-\N$ is a closed domain in $\BHH$ with $\PI \Omega =\Gamma$. Now, we will show that $\Omega$ is a mean convex domain in $\BHH$ and finish the Step 2a. In order to see this, first note that the smooth parts of $\partial \Omega$ are minimal surfaces, so they are already mean convex. So, the only issue is the nonsmooth parts in $\partial \Omega^+$. We need to see that the dihedral angle at the nonsmooth parts are less than $\pi$. We can see this as follows: In the $\partial \wt{\N}$, the nonsmooth parts are coming from intersection of area minimizing surfaces bounding rectangles. These surfaces continue inside $\wt{\N}$ in each intersection, so they have dihedral angles less than $\pi$. For the intersection of infinite side barriers and Scherk barriers, we have infinite arcs $\{\sigma_{i_j}\}$. Notice that in our construction for Scherk barriers, the $\beta$-curves of exact polygons stays either in $\U_i$ lens regions, or $\V_j$ polygons in $\wh{\Delta}^+-\Delta^+$. In both cases, Scherk surfaces smoothly continue inside the $\wh{\U}_i$ or $\W_j$, meaning they smoothly continue inside $\N$. Again, this implies that the dihedral angles along these intersection curves are less than $\pi$. This shows that $\Omega$ is a mean convex domain in $\BHH$.

Now, we will go over remaining trivial cases: If $\Gamma^\pm_g=\emptyset$, then we define $\N^\pm$ without the Scherk barriers $\Y_i^\pm$. Similarly, if either $m=1$ or $n=1$, we can define the barrier $\N^\pm$ without the Scherk barriers $\Y_i^\pm$. Hence in these cases, $\N^\pm=\bigcup\wh{\U}^\pm_k\bigcup\W_j^\pm$. 

If $\Gamma^\pm_c=\emptyset$, then we define $\N^\pm$ without the Scherk barriers $\W_i^\pm$. On the other hand, if $\Gamma^\pm=\Gamma^\pm_c$ is just one point, then define $\N^\pm=\wh{\U}^\pm$.

Finally, if $\Gamma^+=\emptyset$, then let $t^+_0=\sup\{t\in\BR\mid (\theta,t)\in\Gamma\}<\infty$ be the highest height of $\Gamma$. Then define $\N^+=\BH^2\times (t_0,+\infty)$. Similarly, if $\Gamma^-=\emptyset$, define $\N^-=\BH^2\times (t_0^-,-\infty)$ where $t_0^-$ is the lowest height of $\Gamma$.  Step 2a follows. \hfill $\Box$

\vspace{.2cm}

In the following step, we will construct our sequence of compact area minimizing surfaces $\{\Sigma_n\}$ in $\Omega$ so that it cannot escape to infinity.

\vspace{.2cm}

\noindent{\bf Step 2b:} The sequence of compact area minimizing surfaces $\{\Sigma_n\}$.

\vspace{.2cm}

Let $B_n$ be the $n$-disk in $\BH^2$ with the center origin $O$. Let $\B_n=B_n\times[-n,n]$ be the solid cylinder with height $2n$ and radius $n$. Then $\partial\B_n$ is the cylinder with caps.

Let $\Psi_n$ be a radial projection from $\SI$ into $\partial\B_n$ which maps the corner circles of $\SI$ into corner circles of the cylinder $\partial\B_n$. Then, define $\gamma_n=\Psi_n(\Gamma)$ which is a collection of Jordan curves in $\partial\B_n$. Notice that $\B_n$ is convex by construction. Then, by solving the Plateau problem in $\B_n$ for $\gamma_n$, we get area minimizing surfaces $\Sigma_n$ in $\B_n$ with $\partial \Sigma_n=\gamma_n$ \cite{Fe}. Notice that since $\B_n$ is a convex domain in $\BHH$, $\Sigma_n$ is smooth, and it is area minimizing not only in $\B_n$, but also in $\BHH$.

Recall that $\Omega=\BHH-\N$ and $\PI\Omega=\Gamma$. Hence, by modifying $\gamma_n$ if necessary, we will assume that $\gamma_n\subset \Omega\cap \partial \B_n$ and $\gamma_n\to\Gamma$. Again, $\Sigma_n$ is the area minimizing surface in $\BHH$ with $\partial\Sigma_n=\gamma_n$.

\vspace{.3cm}

\noindent {\em The Limit of the sequence:} We claim that $\Sigma_n\subset \Omega$. Indeed, our proof works not only for area minimizing surfaces, but also minimal surfaces. In other words, we will show that if $\gamma_n\subset \Omega$ and $S_n$ is a {\em minimal surface} with $\partial S_n=\gamma_n$, then $S_n\subset \Omega$.

By \cite{Co}, we already know that $\Sigma_n\cap \wt{\N}=\emptyset$ as we can foliate the domain separated by finite tall rectangles by minimal surfaces. Hence, if we show that $\Sigma_n\cap \N^\pm=\emptyset$, we are done. Recall that $\N^+=\bigcup_i \W_i \bigcup_j \Y_j \bigcup_k \wh{\U}_k$. 

First, we consider infinite side barriers $\wh{\U}^\pm$. In particular, $\wh{\U}^+_j=\bigcup \wh{\U}^+_{jt}$ where $\wh{\U}^+_{jt}$ is bounded by an infinite rectangle. We claim that if  $\gamma_n\cap \wh{\U}^+_{jt}=\emptyset$, then $\Sigma_n\cap \wh{\U}^+_{jt}=\emptyset$. Recall that $\wh{\U}^+_{jt}$ is the domain in $\BHH$ separated by the area minimizing surface $\T^+_{jt}$ where $\PI\T^+_{jt}$ is an infinite rectangle $\R^+_{jt}$. Then, $\wh{\U}^+_{jt}$ is foliated by the vertical translations of $\T^+_{jt}$. Since $\partial \Sigma_n\cap \wh{\U}^+_{jt}=\emptyset$, if $\Sigma_n\cap \wh{\U}^+_{jt}\neq \emptyset$, then the last point of touch with the minimal foliation will give a contradiction by the maximum principle. This implies $\Sigma_n\cap \wh{\U}^+_{jt}=\emptyset$, and hence $\Sigma_n\cap \wh{\U}^+_j=\emptyset$.

Now, we consider Scherk barriers $\W_i$ and $\Y_j$. Recall that $\W^+_i=\bigcup \W^+_{ik}$ where $\partial \W^+_{ik}\subset \wh{S}^t_{ik}\bigcup \partial \wh{\U}^+$. Recall that $\wh{S}^t_{ik}$ is a subsurface of the Scherk graph $S^t_{ik}$, and the vertical translations of $\wh{S}^t_{ik}$ foliates $\W^+_{ik}$. By above, we already know that if $\gamma_n\cap \wh{\U}^+=\emptyset$, then $\Sigma_n\cap \wh{\U}^+=\emptyset$. Hence, if $\Sigma_n\cap \W^+_{ik}\neq \emptyset$, then the last point of touch with the minimal foliation must be in the interior of $\Sigma_n$. This gives a contradiction by the maximum principle as before. This shows that $\Sigma_n\cap \W^+_i=\emptyset$. Same foliation argument also works for $\Y^+_j$, and we have $\Sigma_n\cap \Y^+_j=\emptyset$. This implies $\Sigma_n\cap\N^+=\emptyset$. Hence, if $\gamma_n\subset \Omega$, then $\Sigma_n\subset \Omega$ as claimed.

In particular, this implies that $\N^\pm$ is indeed a barrier for the sequence $\{\Sigma_n\}$. By a standard compactness theorem of geometric measure theory \cite{Fe}, the sequence of area minimizing surfaces $\{\Sigma_n\}$ has a subsequence converging  in compact sets $\B_n$. Hence, by using the diagonal sequence argument, we obtain a limit area minimizing surface $\Sigma$ in $\BHH$. Furthermore, $\Sigma\subset \Omega$ as $\Sigma_n\subset \Omega$ for any $n$. This shows that $\PI\Sigma\subset \PI\Omega$. As in \cite[Theorem 2.14]{Co}, for the cylinder part $\wt{\Gamma}$, we can use tall rectangles as barriers. This proves that $\PI\Sigma\subset \Gamma$. Since $\partial\Sigma_n=\gamma_n\to\Gamma$, by using the linking argument in \cite{Co}, we conclude that $\PI\Sigma=\Gamma$. Step 2 follows. 
\end{pf}

\begin{rmk} [Nonoverlapping vs. Fillability] Notice that the proof of Step 1c also shows that if a curve $\Gamma$ is overlapping at the corner, it cannot bound a \underline{minimal surface} either. In particular, if $\Gamma$ contains an interval $(\theta_1,\theta_2)\times\{+\infty\}$, then by the proof of the lemma, any minimal surface $S$ with $\PI S=\Gamma$ must intersect the area minimizing surface $\T_t$ bounding a tall rectangle $\R_t$ which is very high. However,  if $\Omega_t$ is the region in $\BHH$ separated by $\T_t$, then we can foliate $\Omega_t$ by minimal surfaces $\{\Sigma_s\}$ \cite[Lemma 2.11]{Co}. Hence, if $S\cap \T_t\neq \emptyset$, then there must be a last point of touch with $S$ and $\{\Sigma_s\}$ which contradicts to maximum principle. We would like to thank the referee for pointing out this observation.
\end{rmk}

\begin{rmk} \label{tamermk} [Wild Curves]
	Notice that in the theorem above, the same proof would work for a special class of wild curves. In particular, if we assume the complement of $\Gamma^\pm_c$ is a union of a suitable countably many intervals in $S^1_\infty\times\{+\infty\}$, then the whole proof goes through even though $\Gamma^\pm_c$ consists of infinitely many points. However, as the Example 2 in Section \ref{wild} shows, the tameness condition for the theorem is crucial, and it is not true for wild curves in general, e.g. when $\Gamma^\pm_c$ is Cantor-like set.
\end{rmk}

\section{Fillable and Non-Fillable Infinite Curves} \label{fillSec}

In this part, we will see that fillability question and strong fillability question are quite different for infinite curves. In particular, for a given infinite tall curve $\Gamma$, while $\Gamma^\pm$ completely determines strong fillability of $\Gamma$, $\Gamma^\pm$ is not very useful to detect whether $\Gamma$ is fillable.


\subsection{$\Gamma^\pm$ and Fillability} \label{fillsec} \

A trivial observation is that for any given collection of disjoint geodesics $\gamma_1\cup...\cup\gamma_n$ in $\BH^2\times\{+\infty\}$, there is a fillable curve $\Gamma$ with $\Gamma^+=\gamma_1\cup...\cup\gamma_n$. In particular, if $\p_i=\gamma_i\times \BR$ is the vertical plane over $\gamma_i$, then $S=\p_1\cup...\cup\p_2$ is a collection of minimal planes, and $\Gamma=\PI S$ would be a fillable curve with $\Gamma^+=\gamma_1\cup...\cup\gamma_n$. So, only $\Gamma^+$ (or only $\Gamma^-$) is not enough to determine if $\Gamma$ is fillable or not.

Furthermore, we will show that knowing both $\Gamma^+$ and $\Gamma^-$ together is not enough to determine whether $\Gamma$ is fillable or not. First, by using the following theorem, we have a very large family of fillable curves.

\begin{thm} \label{fillthm} Let $\Gamma$ be a curve with $n$ components in $\SI$, i.e. $\Gamma=\Gamma_1\cup\Gamma_2\cup..\Gamma_n$. If each $\Gamma_i$ is tall, and fat at infinity, then $\Gamma$ is fillable.
\end{thm}

\begin{pf} By Theorem \ref{main}, there exists an area minimizing surface $\Sigma_i$ with $\PI\Sigma_i=\Gamma_i$. Let $\Sigma=\Sigma_1\cup...\cup\Sigma_n$. Therefore, if $\Sigma_i\cap\Sigma_j=\emptyset$, then $\Sigma$ would be an embedded minimal surface in $\BHH$, and we are done.

By assumption $\Gamma_i\cap\Gamma_j=\emptyset$. Since $\SI$ is topologically a sphere, $\Gamma_i$ and $\Gamma_j$ bounds disjoint regions in $\SI$, i.e. $\Gamma_i=\partial \Omega_i$ and $\Omega_i\cap\Omega_j=\emptyset$. Then by \cite[Lemma 2.19]{Co}, the area minimizing surfaces bounding such disjoint curves are disjoint, i.e. $\Sigma_i\cap\Sigma_j=\emptyset$.

This shows that $\Sigma$ is a collection of disjoint area minimizing surfaces. Therefore, $\Sigma$ is an embedded minimal surface in $\BHH$ with $\PI\Sigma=\Gamma$, and $\Gamma$ is fillable. The proof follows.
\end{pf}

\begin{rmk} \label{notstrong} [Fillable but not Strongly Fillable Curves] Note that the theorem above does not show that such a $\Gamma$ is strongly fillable. This is because the union of area minimizing surfaces may not be area minimizing. For example, if $\Gamma^+$ is skinny at infinity, the surface $\Sigma=\Sigma_1\cup...\cup\Sigma_n$ is not area minimizing by Theorem \ref{main}, even though each $\Sigma_i$ is area minimizing. However, it is still a minimal surface, and hence $\Gamma$ is fillable.

By using this idea, it is easy to construct many examples of fillable, but not strongly fillable curves. In particular, by choosing a collection of Jordan curves $\Gamma=\Gamma_1\cup\Gamma_2..\cup\Gamma_n$ in $\SI$ where $\Gamma$ is skinny at infinity, then $\Gamma$ is fillable, but not strongly fillable by Theorem \ref{main} and Theorem \ref{fillthm}.
\end{rmk}

By using this theorem and the infinite rectangles (Lemma \ref{infrec}), we show that for any given collection of geodesics $\Gamma^+$ and $\Gamma^-$, it is possible to construct a fillable curve $\Gamma$.

\begin{cor} \label{fillcor} Let $\gamma_1^+\cup..\cup\gamma_n^+$ and $\gamma_1^-\cup..\cup\gamma_m^-$ be given collections of disjoint geodesics in $\BH^2\times\{+\infty\}$ and $\BH^2\times\{-\infty\}$ respectively. Then, there exists a fillable curve $\Gamma$ with $\Gamma^+=\gamma_1^+\cup..\cup\gamma_n^+$ and $\Gamma^-=\gamma_1^-\cup..\cup\gamma_m^-$.
\end{cor}

\begin{pf} Let $\PI \gamma_i^+=\{p_i^+,q_i^+\}$, and $\alpha_i^+$ be the shorter arc in $S^1_\infty$ with the endpoints $\{p_i^+,q_i^+\}$. Let $t_1,t_2,..,t_n$ be real numbers such that if $\alpha_i^+\subset\alpha^+_j$, then $t_i>t_j$. Then, let $\wh{\alpha}_i^+=\alpha_i^+\times\{t_i\}$ be an arc in $S^1_\infty\times\{t_i\}$. Then let $\R_i^+$ be the infinite rectangle with $\PI \R_i^+\supset \gamma_i^+\cup\wh{\alpha}^+_i$. Similarly, define $\R_j^-$. As both $\{\gamma^+_i\}$ and $\{\gamma^-_j\}$ are pairwise disjoint family of geodesic arcs, and both $\{\wh{\alpha}_i^+\}$ and $\{\wh{\alpha}_j^-\}$ are arcs in different levels $S^1_\infty\times \{t_i\}$, then $\Gamma=\bigcup_i\R_i^+\bigcup_j\R_j^-$ is a union of pairwise disjoint infinite rectangles (See Figure \ref{fig2}-left).

\begin{figure}[b]
	\begin{center}
		$\begin{array}{c@{\hspace{.4in}}c}

		\relabelbox  {\epsfysize=2in \epsfbox{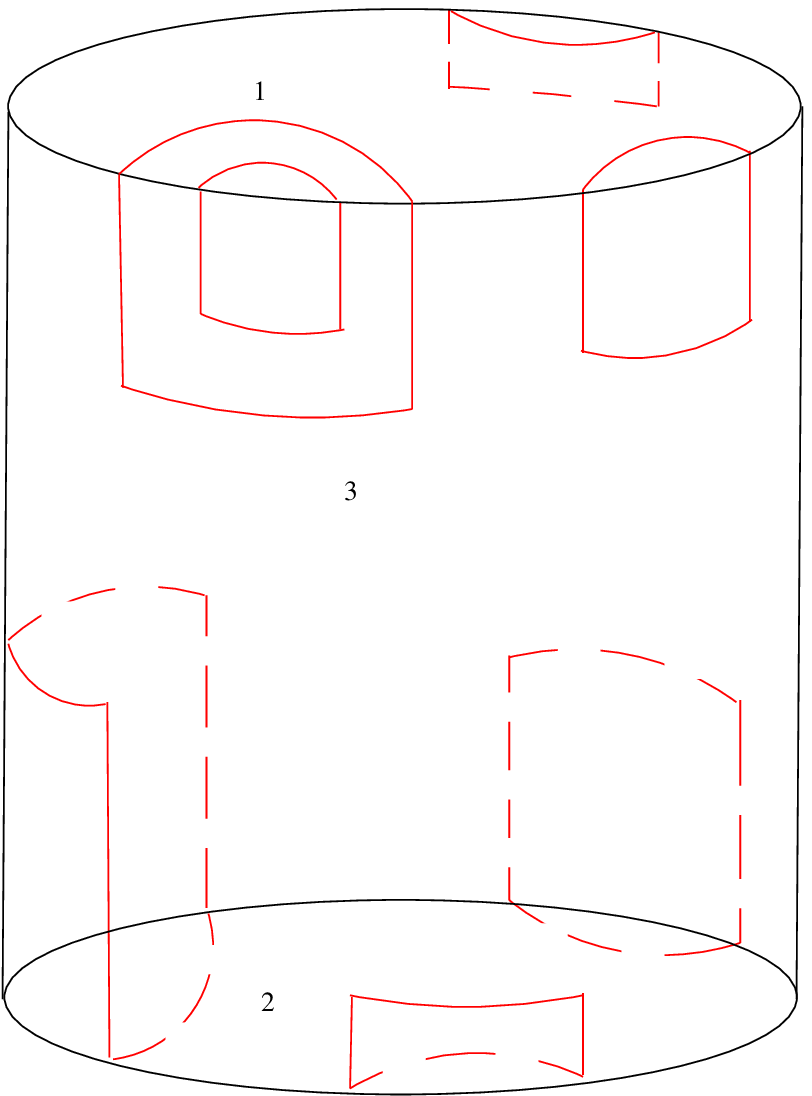}} \relabel{1}{\footnotesize $\Gamma^+$} \relabel{2}{\footnotesize $\Gamma^-$} \relabel{3}{\footnotesize $\R_2^+$} \endrelabelbox &
		
		\relabelbox  {\epsfysize=2in \epsfbox{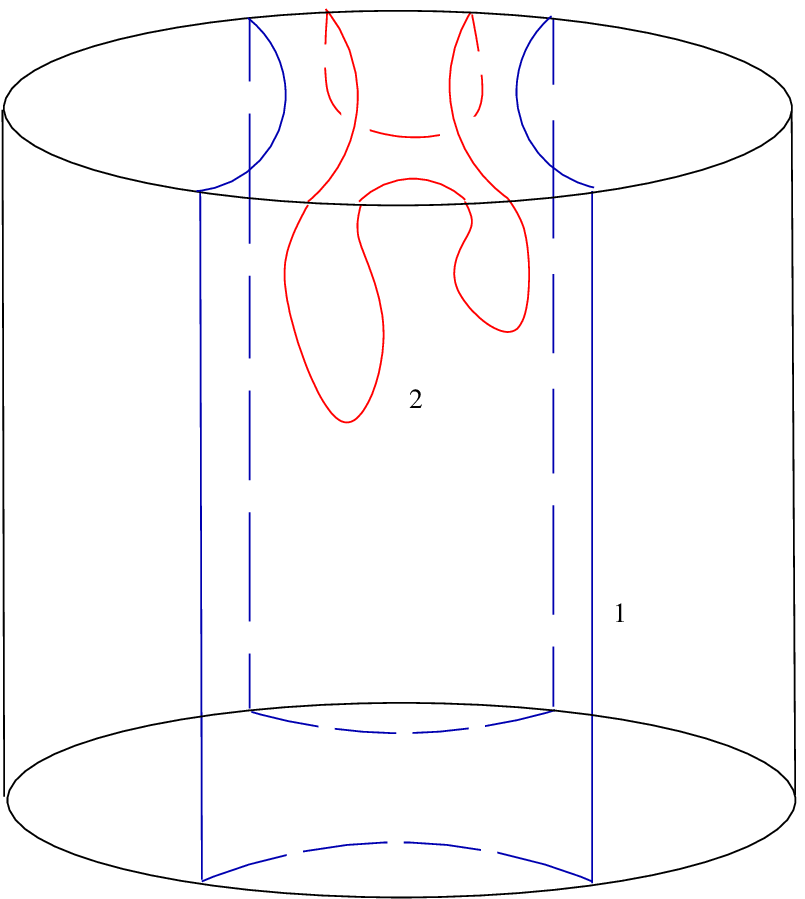}}
		\relabel{1}{\footnotesize $\xi$} \relabel{2}{\footnotesize $\Gamma$ }     \endrelabelbox \\
		\end{array}$
		
	\end{center}
	
	\caption{ \label{fig2} \footnotesize In the figure left, given $\Gamma^+=\gamma_1^+\cup..\cup\gamma^+_4$ and $\Gamma^-=\gamma^-_1\cup..\cup\gamma^-_3$, we construct fillable $\Gamma$ with infinite rectangles $\R^\pm_i$ for each $\gamma^\pm_i$. In the figure right, $\xi$ is a Scherk curve, and $\Gamma$ is trapped by $\xi$.}
\end{figure}

Let $\T_i^+$ be the unique area minimizing surface $\R_i^+$ bounds, and $\T_j^-$ be the unique area minimizing surface $\R_j^-$ bounds. Then, $\s=\bigcup_i\T_i^+\bigcup_j\T_j^-$ is a union of pairwise disjoint minimal surfaces by Theorem \ref{fillthm}. Hence, $\s$ is a complete embedded minimal surface with $\PI \s = \Gamma$. Hence, $\Gamma$ is fillable, and by construction $\Gamma^+=\gamma_1^+\cup..\cup\gamma_n^+$ and $\Gamma^-=\gamma_1^-\cup..\cup\gamma_m^-$. The proof follows.
\end{pf}

\begin{rmk} \label{fillrmk} Notice again that the fillable $\Gamma$ constructed above may not be strongly fillable by Remark \ref{notstrong}. Furthermore, even if $\Gamma$ is fat at infinity, we cannot conclude that
the surface $\s$ constructed in the corollary above is area minimizing. This is because the area minimizing surface $\Sigma$ with $\PI \Sigma=\Gamma$ may be connected, or some other surface in $\BHH$.  In any case, $\s$ is a minimal surface in $\BHH$ with $\PI \s=\Gamma$.
\end{rmk}


\vspace{.2cm}



\subsection{Examples of Non-fillable Curves} \

\vspace{.2cm}

First, of course, if $\Gamma^\pm$ is not a collection of geodesics, then by Lemma \ref{geod}, $\Gamma$ is not fillable. Therefore, we only consider the nontrivial case that $\Gamma^\pm$ is a collection of geodesics.

Even though the previous section shows that $\Gamma^\pm$ fails to detect fillability of $\Gamma$, in some cases, with some conditions on $\wt{\Gamma}$, it shows that $\Gamma$ is non-fillable. In this part, we will give a family of examples of non-fillable infinite curves in $\SI$.

Let $\xi$ be a Scherk curve in $\SI$ (Remark \ref{scherkrem}). Let $\Delta_\xi=\Delta^+_\xi=\Delta^-_\xi$ be  ideal polygon in $\caps$ induced by $\xi^\pm$ (Section \ref{fatsec}). Let $\wh{\Delta}^+_\xi$ be the component of $\PI-\xi$ containing $\Delta^+_\xi$. Similarly, define $\wh{\Delta}^-_\xi$.

\begin{defn} [Trapped Curves] \label{trap} Let $\Gamma$ be an infinite Jordan curve in $\SI$. Let $\Delta^+_\Gamma$ be the component of $\SI-\Gamma$ containing the ideal polygon in $\BH^2\times\{+\infty\}$ induced by $\Gamma^+$. If there exists a Scherk curve $\xi$  such that $\overline{\Delta^+_\Gamma}\subset \wh{\Delta}^+_\xi$, then we will call $\Gamma$ {\em trapped by} $\xi$ with notation $\Gamma\prec\xi$. Similarly, extend the definition to the curves by replacing all $+$ signs with $-$ signs in the corresponding places. See Figure \ref{fig2}-right.
\end{defn}


\begin{thm} \label{trap1} Let $\Gamma$ be an infinite curve trapped by a Scherk curve $\xi$. Then, $\Gamma$ is not fillable.
\end{thm}

\begin{pf} The proof is a straightforward application of the maximum principle. Assume that there exists a minimal surface $\Sigma$ with $\PI\Sigma=\Gamma$. Let $S_0$ be a Scherk graph with $\PI S_0=\xi$. Parametrize all Scherk graphs bounding $\xi$ such that for $t\in \BR$, $S_t=S_0+t$ vertical translation of $S_0$ by $t$. Then, for any $t$, $\PI S_t=\xi$, and the family $\{S_t\}$ foliates the convex region $\Delta_\xi\times\BR$ in $\BHH$.

Without loss of generality, assume $\Delta^+_\Gamma\subset \wh{\Delta}^+_\xi$. Notice that $S_{-t}\to \xi^+\times\BR$ as $t\to \infty$ and $\Gamma\cap\xi=\emptyset$. Then, for sufficiently large $t_o>0$, $S_{-t_o}\cap\Sigma=\emptyset$. As $\Delta^+_\Gamma\cap\Delta_\xi\neq \emptyset$, $\Sigma\cap (\Delta_\xi\times\BR)\neq \emptyset$. Then, let $t_1=\inf\{t\mid S_t\cap\Sigma\neq \emptyset\}$. Then, being the first point of touch, $S_{t_1}$ intersects $\Sigma$ tangentially with lying in one side. This contradicts to the maximum principle. Similar arguments work for $\Delta^-_\Gamma\subset \wh{\Delta}^-_\xi$ case, too. The proof follows.
\end{pf}

\section{Final Remarks}

\subsection{Con-Fillable Curves} \label{confill} \

\vspace{.2cm}

As indicated in Section \ref{fillsec}, the behavior of $\Gamma^\pm$ is highly inadequate to detect fillability of an infinite curve. On the other hand, the main result, Theorem \ref{main}, shows that $\Gamma^\pm$ completely determines the strong fillability of an infinite curve in $\SI$. Notice that the minimal surfaces constructed in Corollary \ref{fillcor} for a given $\Gamma^\pm$ are collection of disjoint minimal surfaces for each component in $\Gamma^\pm$. Hence, the following question becomes very interesting:
{\em If we restrict connectedness on $\Gamma$, or the filling minimal surface, does $\Gamma^\pm$ still plays a crucial role to determine fillability of $\Gamma$?}

This question suggests the following notion between fillability and strong fillability. We will call a curve $\Gamma$ in $\SI$ {\em con-fillable}, if the filling minimal surface $\Sigma$ is connected. i.e. $\Gamma$ is con-fillable if there exists a connected, complete, embedded minimal surface $\Sigma$ in $\BHH$ with $\PI \Sigma=\Gamma$.

\vspace{.2cm}

\begin{question} [Con-fillability] \label{conq} Which infinite curves in $\SI$ bound a connected, complete, embedded minimal surface in $\BHH$ ?
\end{question}

On the other hand, we can state a simpler version of this question:

\begin{question}  \label{conq2} Let $\Gamma$ be an infinite Jordan curve in $\SI$. Is $\Gamma$ fillable?
\end{question}

Notice that in the other versions of the problem, we assume $\Gamma$ to be a finite collection of disjoint Jordan curves in $\SI$. If we assume $\Gamma$ to be one Jordan curve, then the filling surface would automatically be connected. Hence, Question \ref{conq2} is just a simpler case of Question \ref{conq}. Also, by changing Question \ref{conq2} slightly, "\textit{Which Jordan curves in $\SI$ bounds a minimal (or least area) plane in $\BHH$?}" is another interesting question where the analogous question for $\BH^3$ was studied by \cite{An,Co2}.

Both of these questions are very interesting as the examples constructed in Theorem \ref{fillcor} do not apply to these cases. Hence, $\Gamma^\pm$ being fat or skinny might play a crucial role to detect con-fillability of the curve $\Gamma$.

On the other hand, con-fillability and strong fillability are very different notions, where one does not include the other one. In particular, we will give two families of examples of curves which are only con-fillable, and only strongly fillable.

\vspace{.2cm}

\noindent {\em Con-fillable, but not strongly fillable curves:} Horizontal catenoids $\s_2$ and minimal $k$-noids $\s_k$ constructed in \cite{MoR, Py} give an important family of examples for con-fillable curves for any $k\geq 2$. Notice that the asymptotic boundary of horizontal catenoids and minimal $k$-noids consists of $k$ infinite "vertical" Jordan curves, i.e. $\Gamma_k=\PI\s_k=\gamma_1\cup...\gamma_k$ where $\gamma_i$ is the asymptotic boundary of a vertical geodesic plane in $\BHH$. Furthermore, by construction, $\Gamma_k$ is skinny at infinity for any $k\geq 2$. This shows $\Gamma_k$ does not bound any area minimizing surface by Theorem \ref{main}. Hence, these are also examples of con-fillable, but not strongly fillable curves.

\vspace{.2cm}

\noindent {\em Strongly fillable, but not con-fillable curves:} For a given two disjoint geodesics $\tau_1$ and $\tau_2$ in $\BH^2$, define $\Gamma_i=\PI(\tau_i\times\BR)$. By \cite{MoR}, there is a constant $\eta_0>0$ such that if $d(\tau_1,\tau_2)<\eta_0$ then $\Gamma=\Gamma_1\cup\Gamma_2$ bounds a horizontal catenoid in $\BHH$. On the other hand, if $d(\tau_1,\tau_2)>\eta_0$, the curve $\Gamma$ is an example of a strongly fillable curve, which is not con-fillable. In order to see this, first notice that $\Sigma=(\tau_1\cup\tau_2)\times \BR$ is an area minimizing surface by Theorem \ref{main}. Now, assume that there is a connected minimal surface $T$ with $\PI T=\Gamma$. By $d(\tau_1,\tau_2)>\eta_0$, there is a Scherk curve $\beta$ in the region between $\Gamma_1$ and $\Gamma_2$ in $\SI$ with $\beta\cap\Gamma=\emptyset$ by \cite{MoR}. By translating the Scherk graph $S$, we can assume that $S\cap T=\emptyset$. However, as $T$ is connected, a vertical translation of $S$ must intersect $T$. Hence,  first point of touch gives a contradiction by maximum principle as in Theorem \ref{trap1}.

\subsection{Finite Curves} \

\vspace{.2cm}

In \cite{Co}, we discussed the asymptotic Plateau problem for finite curves, and give a fairly complete classification for strongly fillable curves. In this paper, we completed this classification by giving a characterization for infinite strongly fillable curves.

While strong fillability question has been finished, fillability question for finite and infinite curves are wide open. In particular, we gave examples of fillable and non-fillable curves in Section \ref{fillSec}. Furthermore, the same question for finite curves is also very delicate. By a simple generalization of Theorem \ref{fillthm}, if $\Gamma$ is a collection of strongly fillable (finite or infinite) curves, then it is fillable. On the other hand, the only known family of finite non-fillable curves in $\SI$ is the curves containing thin tails by Lemma \ref{thin1}. However, there are many curves $\Gamma$ with $h(\Gamma)<\pi$ containing no thin tails. In \cite{Co,KM}, some families of fillable examples, namely \textit{butterfly curves}, have been constructed. However, it is still wide open question that which curves are fillable among such curves?

Similarly, con-fillability question is also wide open for finite curves, too. \textit{Which curves in $\SI$ bounds a connected minimal surface?} As discussed in previous section, this question is also very different from the fillability, and the strong fillability questions. In \cite{FMMR}, the authors studied a special case of this problem, namely for minimal annuli. By using the ideas for not con-fillable examples above, it can be showed that if the components of a finite curve are "horizontally apart", then it is not con-fillable.

In particular, if $\Gamma=\gamma_1\cup..\cup\gamma_n$ is a finite curve where $\Gamma$ stays in one side of a Scherk curve, then $\Gamma$ can not con-fillable. Similarly, if a finite curve is "vertically $\pi$-apart", then it is not con-fillable. In other words, if $\Gamma=\gamma_1\cup..\cup\gamma_n$ is a finite curve where some $\gamma_i$ is in $S^1_\infty\times(-\infty, c)$ and all other  $\gamma_j$ is in $S^1_\infty\times (c+\pi, \infty)$, then $\Gamma$ cannot be con-fillable. This is because for the pair of horizontal circles $\alpha_\e=S^1_\infty\times \{c+\e, c+\pi-\e\}$ in $\SI$, there is a catenoid $\C_\e$ with $\PI\C_\e=\alpha_\e$ such that $S\cap\C_\e=\emptyset$ for a given minimal surface $S$ with $\PI S=\Gamma$. In particular, as $\e\to 0$, $\C_\e\to \SI$ by \cite{NSST}. This means if there is a connected minimal surface $S$ bounding $\gamma_1$ and other $\gamma_i$'s, it must go through the neck inside catenoid  $\C_\e$ as $\C_\e$ is very close to $S^1_\infty\times [c+\e, c+\pi-\e]$.  Then, by using horizontal hyperbolic translations $\varphi_t$ of $\C_\e$ towards $S$, we can get a contradiction with maximum principle at the first point of contact of $\varphi_t(\C_\e)$ and $S$.

Let $\Gamma$ be {\em semi-infinite} if only one of $\Gamma^+$ or $\Gamma^-$ is nonempty. Fillability, and con-fillability questions for a semi-infinite curve $\Gamma$ might be detected by the behavior at infinity. For example, "\textit{Are there any fillable semi-infinite Jordan curve in $\SI$ which is skinny at infinity?}" seems an interesting question, and easier case to study for fillability question. By considering Theorem \ref{trap1}, if there exists such a curve, it is a curve not trapped by a Scherk curve. While $\Gamma^\pm$ hardly detects fillability for infinite curves by section \ref{fillsec}, it might be the key property for semi-infinite curves.

\subsection{Wild Curves} \label{wild} \

\vspace{.2cm}

In this part, we give two examples of wild curves, which shows that the tameness condition are essential for the results of the paper. We would like to thank the anonymous referee for these examples.\\

\noindent {\bf Example 1:} {\em There are strongly fillable wild curves with $\Gamma^\pm_g$ contains infinitely many geodesics.}

\begin{figure}[b]
	\begin{center}
		$\begin{array}{c@{\hspace{.4in}}c}

		\relabelbox  {\epsfysize=2in \epsfbox{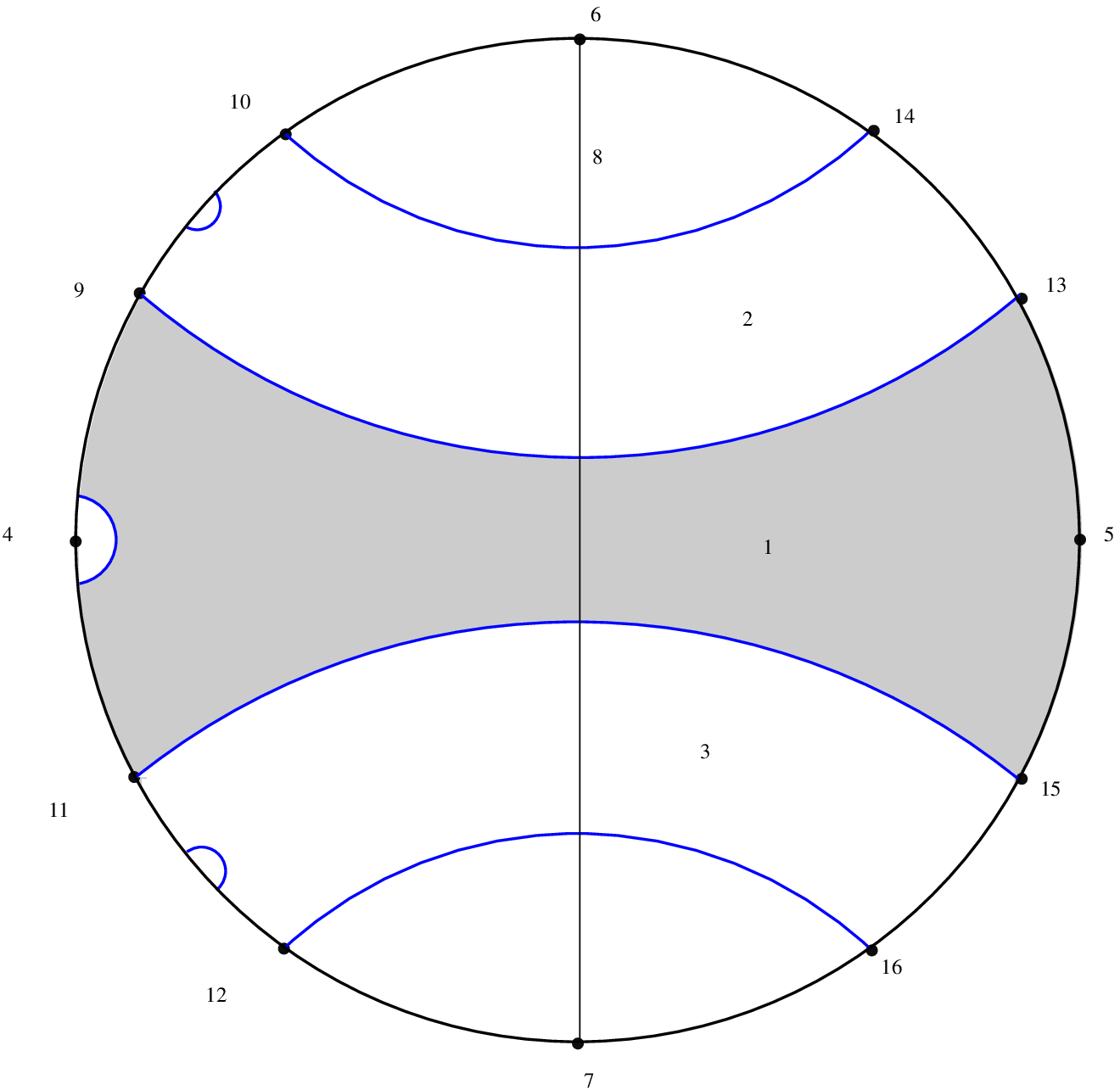}} 
		\relabel{1}{\tiny $\D_0$} \relabel{2}{\tiny $\D_1$} \relabel{3}{\tiny $\D_{-1}$} 
		\relabel{4}{\tiny $p_o$} \relabel{5}{\tiny $q_o$} \relabel{6}{\tiny $a_+$}
		\relabel{7}{\tiny $a_-$} \relabel{8}{\tiny $\gamma$} \relabel{9}{\tiny $p_1$}
		\relabel{10}{\tiny $p_2$} \relabel{11}{\tiny $p_{-1}$} \relabel{12}{\tiny $p_{-2}$}
		\relabel{13}{\tiny $q_1$} \relabel{14}{\tiny $q_2$} \relabel{15}{\tiny $q_{-1}$}
		\relabel{16}{\tiny $q_{-2}$}
		\endrelabelbox &
		
		\relabelbox  {\epsfysize=2in \epsfbox{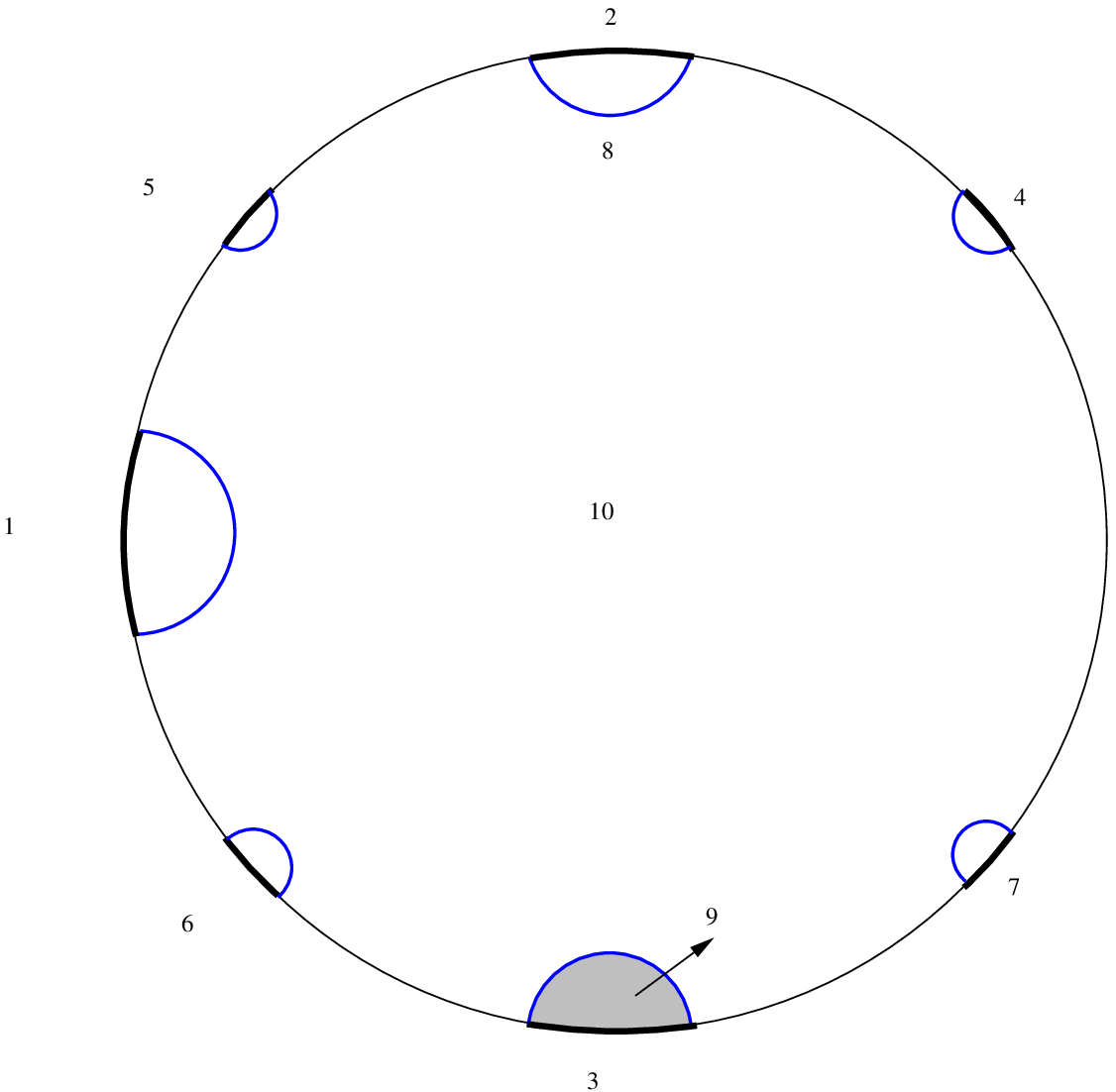}} 
		\relabel{1}{\tiny $I_{1,1}$} \relabel{2}{\tiny $I_{2,1}$} \relabel{3}{\tiny $I_{2,2}$} 
		\relabel{4}{\tiny $I_{3,1}$} \relabel{5}{\tiny $I_{3,2}$} \relabel{6}{\tiny $I_{3,3}$}
		\relabel{7}{\tiny $I_{3,4}$} \relabel{8}{\tiny $\gamma_{2,1}$} \relabel{9}{\tiny $D_{2,2}$}
		\relabel{10}{\scriptsize $\Omega_3$} 		 \endrelabelbox \\
		\end{array}$
		
	\end{center}
	
	\caption{ \label{wildfig} \footnotesize In the figure left, construction of a strongly fillable wild curve $\Gamma$ with $\Gamma^+_g$ has infinitely many geodesics is given. In the figure right, inductive construction of nonfillable wild curve which is nonoverlapping at the corner is shown.}
\end{figure}

Let $p_0\neq q_0\in \PI \BH^2$, and $T$ be a hyperbolic translation along a geodesic $\gamma$ orthogonal to the geodesic $(p_0,q_0)$ where $l$ is the translation length of $T$ along $\gamma$. For $t\in \BR$, we define  $p_t=T^t(p_0)$ and $q_t=T^t(q_0)$ where $T^t$ denotes the translation length $tl$ along $\gamma$. For $t>0$ small, and $l$ large, we consider the domain $\D_0$ bounded by the geodesics $(p_1,q_1), (p_{-1}q_{-1}), (p_{-t},p_{t})$ and the arcs $\overline{p_t p_1}, \overline{q_{-1} q_1}, \overline{p_{-1}, p}_{-t}$  in $\PI \BH^2$ (See Figure \ref{wildfig}-left). If $t$ is small, and $l$ is large, there is a solution $u$ of the minimal surface equation on $\D_0$ with value $+\infty$ on the geodesics, and $0$ on the arcs in $\PI \BH^2$ by \cite{MRR}. We fix the parameter $t$ and $l$. 

Let $\D_k=T^{2k}(\D_0)$, and let $\Omega=\bigcup_{k\in\BZ}\D_k$. Let $a_\pm=\lim_{n\to\pm\infty} T^n(p_0)=\lim_{n\to\pm\infty} T^n(q_0)$.
There are two arcs $\alpha_1$ and $\alpha_2$ in $\PI \BH^2$ joining $a_-$ and $a_+$ with $p_0\in\alpha_1$ and $q_0\in \alpha_2$. Let $f$ be a function defined on $\alpha_1$ and $\alpha_2$ such that $\lim_{p\to a_\pm} f(p)= +\infty$. There is indeed a solution $v$ to the minimal surface equation on $\Omega$ such that $v=+\infty$ on any geodesics $(p_{2k-t},p_{2k+t})$, and $v=f$ on any arcs $\overline{p_{2k+t} p}_{2k+2-t}$ and $\alpha_2$ in $\PI \BH^2$. The existence of such $v$ comes from the fact that $u\circ T^{-2k} +c$ can be used as a barrier from above for the domain $\D_k$. Hence, the graph of $v$ is an area minimizing surface $\Sigma$ with $\PI \Sigma=\Gamma$ such that $\Gamma^+$ is the union of the geodesics $(p_{2k-t}, p_{2k+t})\times\{+\infty\}$ and $(a_\pm,+\infty)$, and $\wt{\Gamma}$ is made of vertical halflines $\{p_{2k\pm t}\}\times[f(p_{2k\pm t},\infty)$, and the graph of $f$ over $\overline{p_{2k+t} p}_{2k+2-t}$ and $\alpha_2$. $\PI \Sigma=\Gamma$ is a wild infinite curve as $\Gamma^+$ has an infinite number of geodesic arcs. This shows that there are strongly fillable wild curves with $\Gamma^\pm_g$ contains infinitely many geodesics.

\vspace{.2cm}

\noindent {\bf Example 2:}	{\em There are nonfillable wild, tall curves which are nonoverlapping at the corner.}
		
The main result states that nonoverlapping tame curves are strongly fillable. This example shows that tameness condition in the main result is crucial. In particular, if we remove tameness condition from $\Gamma$ in the main result, the fillability of $\Gamma$ fails in general. See also Remark \ref{tamermk}.

We will define a curve $\Gamma$ in $\SI$ where $\Gamma^+_c$ contains infinitely many points, but no interval. In particular, $\Gamma^+_c$ will be a Cantor-like set, and the construction will proceed by induction. We will construct a family of disjoint open arcs $I_{n,k}$ in $\PI \BH^2$ where $n\geq 1$ and $k\in\{1,2,...,2^n-1\}$, and an increasing family of compact sets $K_n$ such that the following properties are satisfied. Let $C_n=\bigcup_{m\leq n}\bigcup_{1\leq k\leq 2^m-1} I_{m,k}$.

\begin{itemize}
	\item If $e^{i\frac{k}{2^{n-1}}\pi} \in C_{n-1}$, then $I_{n,k}=\emptyset$. If $e^{i\frac{k}{2^{n-1}}\pi} \notin C_{n-1}$, then $I_{n,k}$ is an arc centered at $e^{i\frac{k}{2^{n-1}}\pi}$ disjoint from $C_{n-1}$ whose endpoints are not of the form $e^{i\frac{p}{m}\pi}$ (See Figure \ref{wildfig}-right).
	
	\item Let $\gamma_{n,k}$ be the geodesic joining the endpoints of $I_{n,k}$, and let $\Omega_n$ be the convex hull of $\PI \BH^2-C_n$ in $\BH^2$. Then, $\bigcup_m K_m=\bigcup_n \Omega_n = \Omega$.
	
	\item There is a solution $u_n$ of the minimal surface equation on $\Omega_n$ whose boundary values are $-\infty$ on any $\gamma_{m,k}$ for any $m\leq n$, and $0$ on $\PI\BH^2-C_n$.
	
	\item $|u_{n+1}-u_n|\leq\frac{1}{2^n}$ on $K_{n+1}$.
	
\end{itemize}

We can see that we can find such a family satisfying above properties as follows: Assume that for $n\leq n_0$, $C_n,K_n, \Omega_n, u_n$ are constructed with the properties above. Then for $n_0+1$, if  $e^{i\frac{k}{2^{n_0}}\pi} \in C_{n_0}$, then choose $I_{n_0+1,k}=\emptyset$, if not, choose $I_{n_0+1,k}$ to be an arc of $\PI \BH^2$ disjoint from $C_{n_0}$. Furthermore, the length of $I_{n_0+1,k}$ can be chosen sufficiently small so that $\gamma_{n_0+1,k}$ is outside of $K_{n_0}$, and $u_{n_0+1}$ exists on $\Omega_{n_0+1}$. Besides, as the lengths of the $I_{n_0+1,k}$'s go to $0$, then $\Omega_{n_0+1}\to \Omega_{n_0}$ and $u_{n_0+1}\to u_{n_0}$ on $\Omega_{n_0}$. This implies we can choose the length of $I_{n_0+1,k}$ such that $|u_{n_0+1}-u_{n_0}|\leq \frac{1}{2^n}$ on some compact set $K_{n_0+1}$ with $K_{n_0}\subset K_{n_0+1}\subset \Omega_{n_0+1}$. Clearly, this can be done so that for any $n,m$, $K_m\subset \Omega_n$. Hence, we get $\bigcup_m K_m=\bigcup_n \Omega_n = \Omega$.

On $\Omega$, $u=\lim u_n$ exists, and it is a solution for minimal surface equation by construction. Since $u\leq u_n$, $u$ takes the value $-\infty$ on any geodesic $\gamma_{n,k}$. Let $C=\bigcup C_n$, and let $E=\PI \BH^2-C$. Then, $\PI \Omega=E$. By construction,  $e^{i\frac{k}{2^{n-1}}\pi} \in C_{n}\subset C$ for any $k,n$. This implies  $e^{i\frac{k}{2^{n-1}}\pi} \notin E$ for any $k,n$, and hence, $E$ contains no interval. A priori, we say nothing about the value of $u$ on $E$ except $u\leq 0$ on $\Omega$.

Let $f$ be defined on $\PI \BH^2$ by $f(p)=\frac{1}{d(p,E)}$ where $d$ is the Euclidean distance along the unit circle. The graph of $f$ is a continuous function on $\PI \BH^2$ with value $\overline{\BR}$. Its graph is a Jordan curve $\Gamma$ in $\SI$. We claim that $\Gamma$ is the desired example.

First, $\Gamma^+_c=E\times\{+\infty\}$, and so $\Gamma$ is nonoverlapping at the corner as $E$ contains no interval. As $\Gamma$ is graph over $\PI \BH^2$, it is also tall. However, we claim that $\Gamma$ is not fillable. Assume that there is a minimal surface $\Sigma$ with $\PI \Sigma=\Gamma$. Since $\Gamma$ is above the closure of the graph $u+t$, $\Sigma$ is also above the graph of $u+t$ for any $t$. By letting $t\to +\infty$, the graphs of $u+t$ will sweep out $\Omega\times\BR$. By maximum principle, $\Sigma\cap \Omega\times\BR=\emptyset$. Hence, $\Sigma$ is contained in $D_{n,k}\times\BR$ where $D_{n,k}$ is the hyperbolic half-space bounded by $\gamma_{n,k}$ and $I_{n,k}$. This means $\Sigma$ has infinitely many components for each $D_{n,k}$ as $\Gamma\cap I_{n,k}\times\BR\neq \emptyset$. 

Let $\Sigma_1$ be the one in $D_{n,k}\times\BR$. By construction, $\PI\Sigma_1=\Gamma\cap I_{1,1}\times\BR$ which is the graph of $f$ over $I_{1,1}$. Let $I_-$ and $I_+$ be two arcs in $\PI \BH^2$ such that $I_-\varsubsetneq I_{1,1}\varsubsetneq I_+$. Let $f_\pm:I_\pm\to \BR$ such that $f<f_-$ on $I_-$, and $f>f_+$ on $I_{1,1}$. Furthermore, we assume $f_\pm=+\infty$ on $\partial I_\pm$. Let $\gamma_\pm$ be geodesic joining endpoints of $I_\pm$ and $D_\pm$ be the hyperbolic half-space bounded by $I_\pm$ and $\gamma_\pm$. Let $\Sigma_\pm$ be the minimal graph over $D_\pm$ with asymptotic boundary the graph of $f_\pm$ over $I_\pm$, and $\gamma_\pm\times\{+\infty\}$. By using the maximum principle, we see that $\Sigma_1$ is above $\Sigma_+$ and below $\Sigma_-$. By letting $I_\pm\to I$ and $f_\pm \to f$, we see that $\Sigma_1$ is actually equal to the minimal graph with asymptotic boundary the graph of $f$ over $I_{1,1}$ and $\gamma_{1,1}\times\{+\infty\}$. This implies $\PI \Sigma$ contains $\gamma_{1,1}\times\{+\infty\}$. This contradicts to the fact that $\PI \Sigma=\Gamma$. 
Hence, $\Gamma$ is an example of a wild curve which is tall, and nonoverlapping at the corner, but not fillable.

\subsection{Asymptotic $H$-Plateau Problem in $\BHH$} \

\vspace{.2cm}

Constant Mean Curvature (CMC) surfaces are natural generalizations of minimal surfaces. In $\BH^3$ and $\BHH$, many analogous questions have been studied in CMC case, too. Like $H\in[0,1)$ for $\BH^3$, $H\in[0,\frac{1}{2})$ is the interesting case for complete embedded $H$-surfaces in $\BHH$.

We call a curve $\Gamma$ in $\SI$ {\em $H$-fillable} if there exists an $H$-surface $\Sigma_H$ in $\BHH$ with $\PI\Sigma_H=\Gamma$. Hence, the following generalization is very natural:

\begin{question}
{\em Which curves in $\SI$ is $H$-fillable for $H\in[0,\frac{1}{2})$?}
\end{question}

Note that this question has been discussed in \cite{NSST}. In particular, by \cite[Theorem 4.1]{NSST}, an $H$-fillable curve $\Gamma$ cannot be finite or semi-infinite. Furthermore, $\wt{\Gamma}$ must be a collection of vertical lines. On the other hand, very different from the $H=0$ case, for $H>0$ case, there are infinite $H$-catenoids, and $H$-paraboloids in $\BHH$ \cite{NSST}. These surfaces make the question very different from the usual asymptotic Plateau problem. Furthermore, because of these a priori properties for $H$-fillable curves mentioned above, and the existence of these $H$-surfaces, it might be a better idea to study this question in the geodesic compactification of $\BHH$ \cite{KM} rather than the product compactification.

Note also that Scherk graphs described in Section \ref{scherksec} were generalized to CMC context, say Scherk $H$-graphs, for $H\in[0,\frac{1}{2})$ by \cite{HRS}. Hence by replacing Scherk graphs with these Scherk $H$-graphs, it might possible to define the fat/skinny at infinity notions for curves in $\SI$, and generalize a version of Theorem \ref{main} to CMC case by following similar ideas. Note that for the nonproper case, with Meeks and Tinaglia, we recently showed the existence of nonproperly embedded $H$-planes in $\BHH$ for any $H\in(0,\frac{1}{2})$ \cite{CMT}.

\end{document}